\documentclass[12pt]{article}
\usepackage{amssymb,amsthm,amsmath,graphicx}

\newtheorem{thm}{Theorem}[section]
\newtheorem{lem}[thm]{Lemma}
\newtheorem{prop}[thm]{Proposition}
\newtheorem{cor}[thm]{Corollary}

\theoremstyle{definition}
\newtheorem{df}[thm]{Definition}
\newtheorem{ex}[thm]{Example}

\numberwithin{equation}{section}

\newcommand{\Homeo}{\operatorname{Homeo}}
\newcommand{\Conv}{\operatorname{Conv}}
\newcommand{\Del}{\operatorname{Del}}
\newcommand{\diam}{\operatorname{diam}}
\newcommand{\ep}{\varepsilon}

\renewcommand{\phi}{\varphi}

\newcommand{\N}{\mathbb{N}}
\newcommand{\Z}{\mathbb{Z}}
\newcommand{\R}{\mathbb{R}}

\title{Orbit equivalence for \\ Cantor minimal $\Z^d$-systems}
\author{Thierry Giordano 
\thanks{Supported in part by a grant from NSERC, Canada} \\
Department of Mathematics and Statistics \\
University of Ottawa \\
585 King Edward, Ottawa, Ontario, Canada K1N 6N5 
\and
Hiroki Matui 
\thanks{Supported in part by a grant 
from the Japan Society for the Promotion of Science} \\
Graduate School of Science \\
Chiba University \\
1-33 Yayoi-cho, Inage-ku, Chiba 263-8522, Japan 
\and
Ian F. Putnam 
\thanks{Supported in part by a grant from NSERC, Canada} \\
Department of Mathematics and Statistics \\
University of Victoria \\
Victoria, B.C., Canada V8W 3P4 
\and 
Christian F. Skau 
\thanks{Supported in part by the Norwegian Research Council} \\
Department of Mathematical Sciences \\
Norwegian University of Science and Technology (NTNU) \\
N-7034 Trondheim, Norway}

\date{}

\begin{document}
\maketitle

\begin{abstract}
We show that every minimal action of any finitely generated abelian group 
on the Cantor set is (topologically) orbit equivalent to an AF relation. 
As a consequence, 
this extends the classification up to orbit equivalence 
of minimal dynamical systems on the Cantor set 
to include AF relations and $\Z^d$-actions. 
\end{abstract}

\section{Introduction}

In this paper we continue the study, 
undertaken in \cite{GPS1}, \cite{GPS3}, \cite{M1}, \cite{M2} and \cite{GMPS1}, 
of the orbit equivalence of minimal dynamical systems on a Cantor set $X$, 
i.e. a compact, totally disconnected, metrizable space 
with no isolated points. 
By a dynamical system, 
we mean to include actions of countable groups 
as well as \'etale equivalence relations 
(we may call an equivalence relation just a relation). 
We recall the definition of \'etale equivalence relations 
(see \cite{R,GPS2} for more information). 
\begin{df}
An equivalence relation $R$ on $X$ 
endowed with a topology $\mathcal{O}$ is said to be \'etale 
if $(R,\mathcal{O})$ is a locally compact Hausdorff $r$-discrete groupoid. 
In other words, $(R,\mathcal{O})$ is a locally compact Hausdorff groupoid, 
and the two canonical projections from $R$ to $X$ are 
open and are local homeomorphisms. 
\end{df}
For an action $\phi$ of a countable discrete group $G$ on $X$ 
by homeomorphisms, 
the orbit relation $R_\phi$ is defined by 
\[
R_\phi=\{(x,\phi^g(x))\in X\times X\mid x\in X, \ g\in G\}, 
\]
and if $\phi$ is free 
(i.e. $\{g\in G\mid\phi^g(x){=}x\}=\{e\}$ for all $x\in X$) 
then $R_\phi$ is \'etale with the topology 
obtained by transferring the product topology on $X\times G$. 
Another rich and tractable class of \'etale equivalence relations is 
the so-called AF equivalence relations, 
which have played an important r\^ole in earlier work 
and will do so in this paper as well. 
Briefly, an \'etale equivalence relation $R$ is AF (or approximately finite) 
if it can be written as an increasing union of compact open subrelations. 
Such relations have a nice presentation 
by means of a combinatorial object called a Bratteli diagram 
(see \cite{GPS1} for more information). 
Recall also that an equivalence relation $R$ on $X$ 
(even without $R$ having any topology itself) is minimal 
if $R[x]=\{y\in X\mid(x,y)\in R\}$ is dense in $X$ for any $x\in X$. 

As already mentioned our primary interest is the notion of orbit equivalence. 
We recall the following definition 
which generalizes the usual definition for group actions. 

\begin{df}\label{definition1}
Let  $X$ and $X'$ be two topological spaces and 
let $R$ and $R'$ be equivalence relations on $X$ and $X'$, respectively. 
We say that $(X,R)$ and $(X',R')$ are orbit equivalent 
if there is a homeomorphism $h:X\to X'$ such that $(h\times h)(R)=R'$. 
\end{df}

In measurable dynamics, the study of orbit equivalence, 
initiated by Dye \cite{D}, was developed by Krieger \cite{K}, 
Ornstein-Weiss \cite{OW} and Connes-Feldman-Weiss \cite{CFW} 
among many others in the amenable case. 
The strategy of their proofs consisted of showing that 
any amenable measurable equivalence relation is orbit equivalent to 
a hyperfinite measurable equivalence relation and 
classifying these ones. 
In the non-singular case, 
the complete invariant of orbit equivalence is an ergodic flow, 
the so-called associated flow (\cite{K}). 
There has also been considerable work done with the same strategy 
in the category of Borel equivalence relations.  
Hyperfiniteness of Borel $\Z^n$-actions was proved by Weiss. 
More generally, 
Jackson-Kechris-Louveau \cite{JKL} showed hyperfiniteness of 
Borel actions of any finitely generated groups with polynomial growth. 
Recently, it was also proved that 
Borel actions of any countable abelian groups are hyperfinite (\cite{GJ}). 
Complete classification of hyperfinite Borel equivalence relations 
up to orbit equivalence was given in \cite{DJK} 
(see also \cite{KM}). 

Our strategy in the topological case is similar: 
\begin{enumerate}
\item[1)] Provide an invariant of (topological) orbit equivalence, 
which is an ordered abelian group whose definition is given in Section 2. 
\item[2)] Show that this invariant is complete 
for minimal AF relations on the Cantor set. 
This classification was obtained in \cite{GPS1} 
(see \cite{P} for a new proof). 
\item[3)] Prove that more general minimal equivalence relations are 
\emph{affable} i.e. orbit equivalent to an AF equivalence relation. 
A fortiori, this extends the classification 
to a larger family of equivalence relations. 
A key technical ingredient for this step is 
the absorption theorem for minimal AF relations, which states that 
a `small' extension of a minimal AF relation is 
orbit equivalent to the AF relation. 
\end{enumerate}
The main contribution of this paper is to establish such a result 
for minimal actions of finitely generated abelian groups on a Cantor set. 

Let $\phi$ be a free minimal action of $\Z^d$ as homeomorphisms 
on the Cantor set $X$ and 
let $R_\phi$ be the associated \'etale equivalence relation. 
Notice that if $\{\mathcal{T}(x)\mid x\in X\}$ is 
a family of tessellations of $\R^d$ with compact cells and 
such that $\mathcal{T}(\phi^n(x))=\mathcal{T}(x)+n$ for $n\in\Z^d$, 
then we can associate to it a finite subrelation $R_\mathcal{T}$ of $R_\phi$ 
by stating that 
\[
(\phi^p(x),\phi^q(x))\in R_\mathcal{T}\text{ if and only if }
p,q \text{ belongs to the same cell of }\mathcal{T}(x). 
\]
Therefore the first step of the proof of the affability of $R_\phi$ is 
to construct a family of nested sequences 
$\{\mathcal{T}_l(x)\mid x\in X\}_{l\geq1}$ 
of tessellations of $\R^d$ such that 
the associated sequence $(R_{\mathcal{T}_l})_{l\geq1}$ 
of finite equivalence relations defines a minimal AF subrelation of $R_\phi$ 
satisfying the assumptions of the absorption theorem (\cite[Theorem 3.2]{M3}). 

For the case $d=2$, this construction (done in \cite{GMPS2}) involved 
a precise control on the geometry of the cells of the tessellations. 
More precisely, given a Delaunay set $P$ 
(i.e. a both uniformly discrete and relatively dense subset) of $\R^2$, 
we had to modify the Voronoi tessellation associated to $P$ 
to ensure that the disjoint cells of the tessellation are well separated. 

For $d\geq3$, three new issues arise: 

a) the first one is the geometry of the tessellations. 
The argument used to modify the Voronoi tessellations in dimension two 
does not work for $d\geq3$. 
For example, given a Delaunay subset $P$ of $\R^3$, 
its Delaunay triangulation may contain slivers. 
These are tetrahedron whose four vertices lie close to a plane and 
whose projection to the plane is a convex quadrilateral with no short edge. 
In \cite{CDEFT}, 
Cheng, Dey, Edelsbrunner, Facello and Teng proved 
the existence of a triangulation whose vertices are $P$ and with no slivers. 
More precisely, they showed that 
there exists an assignment of weights so that 
the weighted Delaunay triangulation (defined in Section 3) 
contains no slivers. 
We generalize their argument to 
triangulations of $\R^d$ for any $d>2$ in Section 3 and 4. 

b) the second one is the combinatorics 
related to the nesting of the tessellations 
which is substantially more complicated 
being a lot less geometrically intuitive. 
This construction of families of nested sequences of tessellations 
is done in Section 6. 

c) the third one is the application of the absorption theorem. 
The first version of this theorem was stated and shown in \cite{GPS2}. 
In \cite{GMPS2}, 
it was generalized to the case that the `small' extension is 
given by a compact relation transverse to the AF equivalence relation. 
This version was sufficient for the study of $\Z^2$-actions, 
but too restricted for successive applications of the theorem. 
Building on the idea of \cite{GMPS2}, 
the absorption theorem was further generalized in \cite{M3} and 
it is this last version 
which is used inductively in Section 7 of this paper. 

The authors would like to thank the referee for many helpful suggestions.

\section{Main results}

Let $\Omega$ be a compact metrizable space and 
let $G$ be a locally compact group. 
A group homomorphism $\phi:G\to\Homeo(\Omega)$ is called an action, 
if $G\times\Omega\ni(g,x)\mapsto\phi^g(x)\in G$ is continuous. 
The action $\phi$ is said to be free, 
if $\phi^g(x)\neq x$ for any $x\in X$ and $g\in G\setminus\{e\}$, 
where $e$ is the neutral element of $G$. 
The action $\phi$ is said to be minimal, 
if $\{\phi^g(x)\in X\mid g\in G\}$ is dense in $X$ for every $x\in X$. 

We denote by $\R^d$ the usual $d$-dimensional Euclidean space. 
For $p\in\R^d$ and $r>0$, 
$B(p,r)$ denotes the open ball of radius $r$ centred at $p$. 

\begin{df}\label{Ctrans}
Let $d\in\N$. 
Let $\phi$ be a free action of $\R^d$ 
on a compact metrizable space $\Omega$. 
We call a closed subset $X\subset\Omega$ a flat Cantor transversal, 
if the following are satisfied. 
\begin{enumerate}
\item $X$ is homeomorphic to a Cantor set. 
\item For any $x\in\Omega$, $\phi^p(x)$ is in $X$ for some $p\in\R^d$. 
\item There exists a positive real number $M>0$ such that 
\[
C=\{\phi^p(x)\mid x\in X,p\in B(0,M)\}
\]
is open in $\Omega$ and 
\[
X\times B(0,M)\ni(x,p)\mapsto\phi^p(x)\in C
\]
is a homeomorphism. 
\item For any $x\in X$ and $r>0$, 
there exists an open neighbourhood $U\subset X$ of $x$ in $X$ such that 
$\{p\in B(0,r)\mid\phi^p(x)\in X\}=\{p\in B(0,r)\mid\phi^p(y)\in X\}$ 
for all $y\in U$. 
\end{enumerate}
\end{df}

We note that the third property means that 
$\Omega$ is locally homeomorphic to the product of $X$ and $\R^d$. 
If $\phi$ is minimal, then (2) follows from (3) automatically. 

When $X\subset\Omega$ is a flat Cantor transversal, 
we define an equivalence relation $R_\phi$ on $X$ by 
\[
R_\phi=\{(x,\phi^p(x))\mid x\in X,\phi^p(x)\in X,p\in\R^d\}. 
\]
If $\phi$ is a minimal action on $\Omega$, 
then $R_\phi$ is a minimal equivalence relation on $X$. 
Indeed, for any $x\in X$ and a non-empty open subset $U\subset X$, 
(3) implies that 
$\{\phi^q(y)\mid y\in U, q\in B(0,M)\}$ is open in $\Omega$ 
and the minimality of $\phi$ implies that 
$\phi^p(x)$ belongs to this open set for some $p\in\R^d$, 
and so $\phi^{p-q}(x)\in U$ for some $q\in B(0,M)$. 
We provide $R_\phi$ with a topology 
whose basis is given by the sets of the form 
$\{(x,\phi^p(x))\mid x\in U\}\cap R_\phi$, 
where $U$ is an open subset of $X$ and $p$ is in $\R^d$. 
It is not hard to see that 
$R_\phi$ is an \'etale equivalence relation with this topology. 

We have two important classes of examples of flat Cantor transversals. 

\begin{ex}
Let $X$ be a Cantor set and 
let $\phi$ be a free minimal action of $\Z^d$ on $X$. 
Let $\Omega$ denote the quotient space of $X\times\R^d$ 
by the equivalence relation 
$\{((x,p),(\phi^n(x),p+n))\mid x\in X,p\in\R^d,n\in\Z^d\}$. 
The topological space $\Omega$ is called a suspension space of $(X,\phi)$. 
There exists a natural $\R^d$-action on $\Omega$ 
induced by the translation $(x,p)\mapsto(x,p+q)$ in $X\times\R^d$. 
We denote this action by $\tilde\phi$. 
Clearly $\tilde\phi$ is free and minimal and 
it is easy to see that $X\times\{0\}$ is a flat Cantor transversal 
for $(\Omega,\tilde\phi)$. 
In addition, by identifying $X\times\{0\}$ with $X$, 
the \'etale equivalence relation $R_{\tilde\phi}$ agrees with 
the \'etale equivalence relation arising from $(X,\phi)$. 
We also remark that 
the $\R^d$-action $\tilde\phi$ is a special case of 
the so-called Mackey action (see Section 4 of \cite{Mac}). 
\end{ex}

\begin{ex}
We would like to explain equivalence relations 
arising from tiling spaces briefly. 
The reader should see \cite{KP} and the references given there 
for more details 
and \cite{BBG,BG} for further developments. 
We follow the notation used in \cite{KP}. 
Let $T$ be a tiling of $\R^d$ satisfying the finite pattern condition and 
let $\Omega_T$ be the continuous hull of $T$. 
The topological space $\Omega_T$ is compact and metrizable. 
There exists a natural $\R^d$-action $\phi$ on $\Omega_T$ 
given by translation. 
Assume further that $T$ is aperiodic and repetitive. 
Then the action $\phi$ is free and minimal. 
Let $\Omega_{\text{punc}}\subset\Omega_T$ be 
as in Definition 5.1 of \cite{KP}. 
It is easy to check that $\Omega_{\text{punc}}$ is a flat Cantor transversal. 
In addition, the \'etale equivalence relation $R_\phi$ defined above 
agrees with $R_{\text{punc}}$ of \cite{KP}. 
\end{ex}

\bigskip

The following is the main theorem of this paper. 

\begin{thm}\label{main}
Let $X\subset\Omega$ be a flat Cantor transversal 
for a free minimal action $\phi$ of $\R^d$ on $\Omega$ and 
let $R_\phi\subset X\times X$ be the minimal equivalence relation 
induced from $(\Omega,\phi)$. 
Then $R_\phi$ is affable. 
\end{thm}

The proof is quite long and 
we defer it until the last section of the paper. 

Let us recall the algebraic invariant 
associated to an equivalence relation on the Cantor set $X$ and 
add a remark on minimal actions of finitely generated abelian groups on $X$. 

Let $R$ be an equivalence relation on a Cantor set $X$. 
We say that a Borel probability measure $\mu$ on $X$ is $R$-invariant, 
if $\mu$ is $\gamma$-invariant for any $\gamma\in\Homeo(X)$ 
satisfying $(x,\gamma(x))\in R$ for all $x\in X$. 
We let $M(X,R)$ denote the set of $R$-invariant probability measures on $X$. 
This is a weak$*$ compact convex set, 
in fact, a non-empty Choquet simplex 
whenever $R$ arises from a free action of an amenable group. 
We say that $(X,R)$ is uniquely ergodic, 
if the set $M(X,R)$ has exactly one element. 
We let $C(X,\Z)$ denote 
the set of all continuous integer-valued functions on $X$. 
It is an abelian group with the operation of pointwise addition. 
The quotient group of $C(X,\Z)$ by the subgroup 
\[
\left\{f\in C(X,\Z)\mid\int_Xf\,d\mu=0\text{ for all }\mu\in M(X,R)\right\}
\]
is denoted by $D_m(X,R)$. 
For a function $f$ in $C(X,\Z)$, 
we denote its class in $D_m(X,R)$ by $[f]$. 
Of course, this is a countable abelian group, 
but it is also given an order structure \cite{GPS1} 
by defining the positive cone $D_m(X,R)^+$ as the set of all $[f]$, 
where $f\geq0$. 
It also has a distinguished positive element $[1]$, 
where $1$ denotes the constant function with value $1$. 
Our invariant is the triple $(D_m(X,R),D_m(X,R)^+,[1])$. 
This ordered abelian group first appeared in \cite{GPS1} 
in the case of minimal $\Z$-actions 
as the quotient of a $K$-group by its subgroup of infinitesimal elements 
and was shown to be a simple dimension group. 
If a homeomorphism $h:X_1\to X_2$ induces an orbit equivalence 
between $R_1$ on $X_1$ and $R_2$ on $X_2$, 
then clearly we have $h_*(M(X_1,R_1))=M(X_2,R_2)$. 
Hence $h$ induces an isomorphism 
between the two triples $(D_m(X_1,R_1),D_m(X_1,R_1)^+,[1])$ 
and $(D_m(X_2,R_2),D_m(X_2,R_2)^+,[1])$. 

Let us add a remark about actions of finitely generated abelian groups. 
Let $G$ be a finitely generated abelian group and 
let $\phi$ be a minimal action of $G$ on a Cantor set $X$. 
It is easy to see that, for any $g\in G$, 
$X_g=\{x\in X\mid\phi^g(x)=x\}$ is a $G$-invariant closed subset. 
By the minimality of $\phi$, 
if $X_g$ is not empty, it must be the whole of $X$. 
Letting $H$ denote the set of all elements $g$ for which $X_g=X$, 
the orbits of $\phi$ can be realized as 
the orbits of a free action of the quotient group $G/H$, 
which is also a finitely generated abelian group. 
Since finite groups cannot act minimally on an infinite space, 
$G/H$ is isomorphic to $\Z^d\oplus K$, 
where $d\geq1$ and $K$ is a finite abelian group. 
In the Appendix, generalizing a result of O. Johansen (\cite{J}), 
we show that any free minimal action of $\Z^d\oplus K$ on a Cantor set $X$ 
is orbit equivalent to a free minimal action of $\Z^d$ 
(see Theorem \ref{app}). 

Then as an immediate consequence of Theorem \ref{main} 
and \cite[Theorem 2.3]{GPS1}, 
we have the following extension of \cite[Theorem 1.6]{GMPS2}. 

\begin{thm}
Let $(X,R)$ and $(X',R')$ be 
two minimal equivalence relations on Cantor sets 
which are either AF relations or 
arise from actions of finitely generated abelian groups 
or tiling spaces on $\R^d$. 
Then they are orbit equivalent if and only if 
\[
(D_m(X,R),D_m(X,R)^+,[1])\cong(D_m(X',R'),D_m(X',R')^+,[1]), 
\]
meaning that there is a group isomorphism between $D_m(X,R)$ and $D_m(X',R')$ 
which is a bijection between positive cones and preserves the class of $1$. 
\end{thm}

As explained in \cite{GMPS1}, 
the range of the invariant $D_m$ for minimal AF relations is precisely 
the collection of simple, acyclic dimension groups 
with no non-trivial infinitesimal elements, 
and exactly the same holds for minimal $\Z$-actions on Cantor sets 
(see \cite{HPS} and \cite{GPS1}). 
It follows from the theorem above that 
every minimal free $\Z^d$-action on a Cantor set is 
orbit equivalent to a $\Z$-action. 
But, we do not have an exact description of the range 
for minimal $\Z^d$-actions, when $d$ is greater than one. 

The following two corollaries are immediate consequences 
of the main theorem and the definitions. 

\begin{cor}
Let $(X,R)$ and $(X',R')$ be two minimal equivalence relations on Cantor sets 
which are either AF relations or 
arise from actions of finitely generated abelian groups 
or tiling spaces on $\R^d$. 
Then they are orbit equivalent if and only if 
there exists a homeomorphism $h:X\to X'$ 
which implements a bijection between the sets $M(X,R)$ and $M(X',R')$. 
\end{cor}

\begin{cor}
Let $(X,R)$ and $(X',R')$ be 
two minimal, uniquely ergodic equivalence relations on Cantor sets 
which are either AF relations or 
arise from actions of finitely generated abelian groups 
or tiling spaces on $\R^d$. 
Suppose that $M(X,R)=\{\mu\}$ and $M(X',R') = \{\mu'\}$. 
Then the two systems are orbit equivalent if and only if 
\[
\{\mu(U)\mid U\text{ is a clopen subset of }X\}
=\{\mu'(U')\mid U'\text{ is a clopen subset of }X'\}. 
\]
\end{cor}

\section{Weighted Delaunay triangulations}

We will be constructing various tessellations of $\R^d$ 
in order to obtain a minimal AF subrelation 
to which the absorption theorem can be applied. 
A similar construction for Cantor minimal dynamical systems 
was first made by Forrest \cite{F}. 
Using Voronoi tessellations, he showed that 
the equivalence relation arising from a free minimal $\Z^d$-action 
has a `large' AF subrelation. 
As mentioned in the introduction, 
for our purposes, however, Voronoi tessellations have a serious drawback. 
We would like to know that 
disjoint cells should be separated in some controlled manner. 
In the case of $\Z^2$-actions, 
we could do this by moving the vertices of the Voronoi tessellation 
to the incentres of the triangles of its dual tessellation. 
But, as noted in \cite{GMPS1}, 
this argument does not work in the situation of $\R^d$ for $d>2$. 

\begin{figure}[h]
\begin{center}
\includegraphics[height=0.2\textheight]{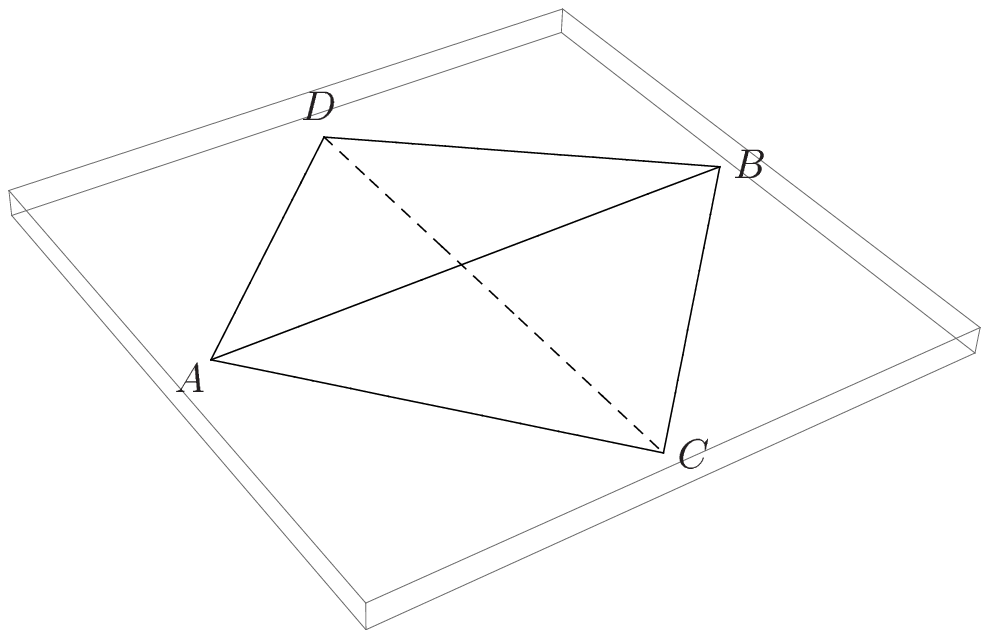}
\caption{sliver in $\R^3$}
\label{fig1}
\end{center}
\end{figure}

The difficulty in the high dimensional space is 
to get rid of badly shaped tetrahedrons 
in the dual of the Voronoi tessellations. 
Such bad tetrahedrons are called slivers. 
For convenience, let us restrict our attention to $\R^3$ for a moment. 
A tetrahedron in $\R^3$ is called a sliver, 
if its four vertices lie close to a plane and 
its projection to the plane is a convex quadrilateral with no short edge 
(see Figure \ref{fig1}). 
In \cite{CDEFT}, using the idea of weighted Delaunay triangulations, 
they studied how to get rid of slivers. 
More precisely, it was shown that 
there exists an assignment of weights so that 
the weighted Delaunay triangulation contains no slivers. 
In this section we generalize their argument to 
triangulations of $\R^d$ for any $d>2$. 
Note that 
the hypothesis used in \cite{CDEFT} is weaker 
than the one needed in our setting, 
and so the result of \cite{CDEFT} is not recovered from our result. 
The reader may refer to \cite{Sch} 
for weighted Delaunay triangulations. 
Actually, in \cite{Sch}, 
the notions of weighted Voronoi and Delaunay tilings were 
unified by Laguerre tilings. 

\bigskip

Let $d(\cdot,\cdot)$ denote the usual Euclidean metric on $\R^d$. 
For any non-empty set $A\subset\R^d$ and $p\in\R^d$, 
we let $d(p,A)=\inf\{d(p,q)\mid q\in A\}$ and $A+p=\{q+p\mid q\in A\}$. 
For a subset $A\subset\R^d$, 
the closed convex hull of $A$ is written by $\Conv(A)$. 
A $k$-dimensional affine subspace of $\R^d$ is 
a translation of a $k$-dimensional vector subspace of $\R^d$. 
A tessellation of $\R^d$ is a collection of compact subsets of $\R^d$ 
which cover $\R^d$ with pairwise disjoint interiors. 
An element of a tessellation is called a cell. 
A tessellation is called a triangulation, 
if every cell of it is a $d$-simplex. 

Let $P\subset\R^d$ be a countable subset. 
For a real number $M$, 
we say that $P$ is $M$-separated 
if $d(p,q)\geq M$ for all $p\neq q\in P$ and 
we say that $P$ is $M$-syndetic 
if $\bigcup_{p\in P}B(p,M)=\R^d$. 
The following lemma is an easy geometric exercise, but 
we include the proof for completeness. 

\begin{lem}\label{sepasyn}
\begin{enumerate}
\item For any $M$-separated set $P\subset\R^d$, $p\in\R^d$ and $R>0$, 
\[
\#(P\cap B(p,R))\leq(2R+M)^d/M^d. 
\]
\item For any $M$-syndetic set $P\subset\R^d$, $p\in\R^d$ and $R>M$, 
\[
\#(P\cap B(p,R))\geq(R-M)^d/M^d. 
\]
\end{enumerate}
\end{lem}
\begin{proof}
Let $V_d$ be the $d$-dimensional volume of the unit ball in $\R^d$. 

(1). Consider the balls $B(q,M/2)$ for all $q\in P\cap B(p,R)$. 
They are mutually disjoint and contained in $B(p,R+M/2)$. 
Hence one has 
\[
\#(P\cap B(p,R))\times(M/2)^dV_d\leq (R+M/2)^dV_d, 
\]
and so $\#(P\cap B(p,R))$ is not greater than $(2R+M)^d/M^d$. 

(2). Consider the balls $B(q,M)$ for all $q\in P\cap B(p,R)$. 
Clearly their union covers $B(p,R-M)$. 
Therefore we get 
\[
\#(P\cap B(p,R))\times M^dV_d\geq (R-M)^dV_d, 
\]
and so $\#(P\cap B(p,R))$ is not less than $(R-M)^d/M^d$. 
\end{proof}

Let $P$ be an $M$-separated and $2M$-syndetic subset of $\R^d$. 
For each $p$ in $P$, let 
\[
V(p)=\{q\in\R^d\mid d(q,p)=d(q,P)\}, 
\]
which is a $d$-polytope with $p$ in its interior. 
The collection $\{V(p)\mid p\in P\}$ is a tessellation of $\R^d$ 
and called the Voronoi tessellation. 

For $k=1,2,\dots,d$, 
we denote by $\Delta_k(P)$ 
the set of all $k+1$ distinct points in $P$ 
which do not lie on a $(k-1)$-dimensional affine subspace. 
In other words, 
$\{p_0,p_1,\dots,p_k\}\subset P$ belongs to $\Delta_k(P)$ 
if and only if 
$\{p_1-p_0,p_2-p_0,\dots,p_k-p_0\}$ is linearly independent. 
For $\tau\in\Delta_d(P)$, 
there exists a unique closed ball whose boundary contains $\tau$. 
We call this closed ball the circumsphere of $\tau$. 
Let $\Del(P)$ be the set of all 
$\tau\in\Delta_d(P)$ such that 
their circumsphere do not intersect with $P\setminus\tau$. 
Generically, $\{\Conv(\tau)\mid \tau\in\Del(P)\}$ 
gives a tessellation of $\R^d$, 
which is called the Delaunay triangulation. 
It is well known that the Delaunay triangulation is 
the dual of the Voronoi tessellation. 

Next, we have to introduce the notion of weighted Delaunay triangulations. 
A pair of $p\in\R^d$ and a non-negative real number $w$ 
is called a weighted point. 
A weighted point $(p,w)$ may be thought as 
a sphere centred at $p$ with radius $\sqrt{w}$. 
We say that 
two weighted points $(p_1,w_1)$ and $(p_2,w_2)$ are orthogonal 
if $d(p_1,p_2)^2=w_1+w_2$. 
Let $p_0,p_1,\cdots,p_k$ be $k+1$ distinct points in $\R^d$ 
which do not lie on a $(k-1)$-dimensional affine subspace, 
where $1\leq k\leq d$. 
Let $H$ be the $k$-dimensional affine subspace containing $p_i$'s. 
Let $w_0,w_1,\cdots,w_k$ be non-negative real numbers such that 
the balls $B(p_i,\sqrt{w_i})$ are mutually disjoint. 
We can see that 
there exists a unique weighted point $(z,u)$ such that 
$z\in H$ and $(z,u)$ is orthogonal to $(p_i,w_i)$ for all $i$ 
in the following way. 
Let $\langle\cdot,\cdot\rangle$ denote the inner product of $\R^d$. 
There exist linearly independent vectors $q_1,q_2,\dots,q_{d-k}$ 
satisfying $\langle p_i{-}p_0,q_j\rangle=0$ for all $i,j$. 
Since $\{p_1{-}p_0,p_2{-}p_0,\dots,p_k{-}p_0\}\cup\{q_1,q_2,\dots,q_{d-k}\}$ 
are linearly independent, 
the system of linear equations 
\[
\langle p_i-p_0,z\rangle
=\frac{1}{2}(\langle p_i,p_i\rangle-w_i-\langle p_0,p_0\rangle+w_0), 
\ i=1,2,\dots,k, 
\]
\[
\langle q_j,z-p_0\rangle=0, 
\ j=1,2,\dots,d-k, 
\]
has the unique solution $z$. 
Put $u=d(p_0,z)^2-w_0$. 
Then $(z,u)$ is orthogonal to $(p_i,w_i)$ for all $i$ and 
$z$ is in $H$. 
Conversely, if $z\in H$ and $(z,u)$ is orthogonal to $(p_i,w_i)$, then 
$z$ must satisfy the system of linear equations above. We call $z$ and $(z,u)$ the orthocentre and the orthosphere of $(p_i,w_i)$, 
respectively. 
Usually the intersection of the three altitudes of a triangle is 
called the orthocentre, 
but we follow the terminology used in \cite{CDEFT}. 

Suppose that 
$P$ is an $M$-separated and $2M$-syndetic subset of $\R^d$ for some $M>0$. 
Put 
\[
{\cal W}=\{w:P\to\R\mid0\leq w(p)\leq (M/3)^2\}. 
\] 
An element $w\in{\cal W}$ is called a weight function on $P$. 

\begin{df}\label{Del}
\begin{enumerate}
\item For $\tau=\{p_0,p_1,\dots,p_k\}\in\Delta_k(P)$, 
we say that $(z,u)$ is the orthosphere of $\tau$ with respect to $w$ 
if $(z,u)$ is the orthosphere of $(p_i,w(p_i))$. 
\item We say that $(z,u)$ is empty in $(P,w)$ 
if $d(z,p)^2-u-w(p)\geq 0$ for all $p\in P$. 
\item We let $\Del(P,w)$ denote the set of all $\tau\in\Delta_d(P)$ 
such that the orthosphere of $\tau$ with respect to $w$ is empty in $(P,w)$. 
\item We let $\mathcal{Z}(P,w)$ be the set of all empty orthospheres 
of some $\tau$ in $\Del(P,w)$. 
\item For $(z,u)\in\mathcal{Z}(P,w)$, 
let $\delta(z,u)$ be the set of all $p\in P$ such that $d(z,p)^2=u+w(p)$. 
\end{enumerate}
\end{df}

Notice that if the weight function $w$ is constantly zero, 
$\Del(P,w)$ agrees with $\Del(P)$. 
The proof of the following can be found in \cite{Sch}. 
We remark that 
the condition R1 in \cite{Sch} is satisfied 
because $P$ is separated and $w$ is bounded, 
and that the condition R2 in \cite{Sch} is satisfied 
because $P$ is syndetic. 

\begin{lem}[\cite{Sch}]\label{Laguerre}
\begin{enumerate}
\item The collection $\{\Conv(\delta(z,u))\mid(z,u)\in\mathcal{Z}(P,w)\}$ is 
a tessellation of $\R^d$. 
\item For each orthosphere $(z,u)\in\mathcal{Z}(P,w)$, 
$\delta(z,u)$ coincides with the extremal points of $\Conv(\delta(z,u))$. 
\item For any orthospheres $(z,u),(z',u')\in\mathcal{Z}(P,w)$, we have 
\[
\Conv(\delta(z,u))\cap\Conv(\delta(z',u'))
=\Conv(\delta(z,u)\cap\delta(z',u')). 
\]
In other words, $\Conv(\delta(z,u))$ and $\Conv(\delta(z',u'))$ meet 
face to face. 
\end{enumerate}
\end{lem}

If $\#\delta(z,u)$ is greater than $d+1$, 
we can partition $\Conv(\delta(z,u))$ 
into $d$-simplices whose vertices are in $\delta(z,u)$ and 
whose interiors are pairwise disjoint. 
In other words, we can find a subset $\mathcal{D}\subset\Del(P,w)$ such that 
$\{\Conv(\tau)\mid\tau\in\mathcal{D}\}$ is a tessellation of $\R^d$. 
We call it the weighted Delaunay triangulation. 

We begin with some basic properties of orthospheres. 
For $\tau\in\Delta_k(P)$, 
we write the $k$-dimensional affine subspace 
which contains $\tau$ by $H(\tau)$. 

\begin{lem}\label{tau'}
Let $w\in{\cal W}$ be a weight function on $P$. 
Suppose that $\tau\in\Delta_l(P)$ is contained in $\tau'\in\Delta_k(P)$, 
where $2\leq l<k$. 
Let $(z,u)$ and $(z',u')$ be the orthospheres 
of $\tau$ and $\tau'$ with respect to $w$, respectively. 
Let $\pi$ be the orthogonal projection from $H(\tau')$ to $H(\tau)$. 
Then we have $\pi(z')=z$ and $d(z',z)^2=u'-u$. 
In particular, $u$ is not greater than $u'$. 
\end{lem}
\begin{proof}
We have $d(z',\pi(z'))^2=d(z',p)^2-d(\pi(z'),p)^2$ for all $p\in\tau$, 
because $z'-\pi(z')$ is orthogonal to $H(\tau)$. 
Combining this with $d(z',p)^2=u'+w(p)$, 
we get $d(\pi(z'),p)^2=u'-d(z',\pi(z'))^2+w(p)$ for all $p\in\tau$. 
We claim that $u'-d(z',\pi(z'))^2$ is non-negative. 
If it were negative, then $d(\pi(z'),p)<\sqrt{w(p)}$ for all $p\in\tau$. 
Thus $\pi(z')\in B(p,\sqrt{w(p)})$. 
Since there are at least two such $p$, 
it does not happen because of $w\in{\cal W}$. 
It follows that $(\pi(z'),u'-d(z',\pi(z'))^2)$ is 
an orthosphere of $(p,w(p))$ for $p\in\tau$. 
By the uniqueness of the orthosphere, 
we get $\pi(z')=z$ and $u'-d(z',\pi(z'))^2=u$, 
which completes the proof. 
\end{proof}

\begin{lem}\label{bddradius}
Let $w\in{\cal W}$ and $\tau=\{p_0,p_1,\dots,p_d\}\in\Del(P,w)$. 
Let $(z,u)$ be the orthosphere of $\tau$ with respect to $w$. 
Then $u$ is less than $4M^2$ and $d(z,p_i)^2<5M^2$. 
In particular, $d(p_i,p_j)$ is less than $2\sqrt{5}M$ 
for any $i,j=0,1,\dots,d$. 
\end{lem}
\begin{proof}
Since $(z,u)$ is empty in $(P,w)$, 
$d(z,p)^2-u-w(p)\geq0$ for all $p\in P$. 
It follows that $B(z,\sqrt{u})$ does not meet $P$, and so 
$\sqrt{u}$ is less than $2M$, because $P$ is $2M$-syndetic. 
For any $i=0,1,\dots,d$, $(z,u)$ is orthogonal to $(p_i,w(p_i))$, 
that is, $d(z,p_i)^2=u+w(p_i)$. 
Therefore $d(z,p_i)^2<5M^2$. 
From this we have $d(p_i,p_j)\leq d(p_i,z)+d(z,p_j)<2\sqrt{5}M$. 
\end{proof}

For each $k=1,2,\dots,d$, 
we let $D_k$ denote the set of all $\tau\in\Delta_k(P)$ 
satisfying that there exist $w\in{\cal W}$ and $\tau'\in\Del(P,w)$ 
such that $\tau\subset\tau'$. 
From the lemma above, 
we can see that, for any $\tau\in D_k$ and $p,q\in\tau$, 
$d(p,q)$ is less than $2\sqrt{5}M$. 

\begin{lem}\label{Dkloc}
For each $k=1,2,\dots,d$, 
$D_k$ is locally derived from $P$ in the following sense. 

For any $p_1,p_2\in P$ satisfying 
\[
(P-p_1)\cap B(0,2\sqrt{5}M)=(P-p_2)\cap B(0,2\sqrt{5}M), 
\]
if $\tau\in D_k$ contains $p_1$, then $\tau-p_1+p_2\in D_k$. 
\end{lem}
\begin{proof}
Suppose that $p_1,p_2\in P$ satisfy 
\[
(P-p_1)\cap B(0,2\sqrt{5}M)=(P-p_2)\cap B(0,2\sqrt{5}M)
\]
and that $\tau\in D_k$ contains $p_1$. 
There exists $w\in{\cal W}$ and $\tau'\in\Del(P,w)$ such that 
$\tau\subset\tau'$. 
From Lemma \ref{bddradius}, 
for any $q\in\tau'$, $d(q,p_1)<2\sqrt{5}M$. 
It follows that $q-p_1+p_2$ belongs to $P$, and so 
$\tau'-p_1+p_2\in\Delta_d(P)$. 
It suffices to show $\tau'-p_1+p_2\in\Del(P,w')$ for some $w'\in{\cal W}$. 
Let $(z,u)$ be the orthosphere of $\tau'$ with respect to $w$. 
As $\tau'\in\Del(P,w)$, 
we have $d(z,p)^2-u-w(p)\geq0$ for all $p\in P$. 
Define $w'\in{\cal W}$ by 
\[
w'(p)=\begin{cases}w(p+p_1-p_2) & p\in B(p_2,2\sqrt{5}M) \\
0 & \text{otherwise}. \end{cases}
\]
for $p\in P$. 
Note that if $p$ is in $P\cap B(p_2,2\sqrt{5}M)$, then 
$p+p_1-p_2$ is in $P\cap B(p_1,2\sqrt{5}M)$. 
Clearly $(z-p_1+p_2,u)$ is the orthosphere of $\tau'-p_1+p_2$ 
with respect to $w'$. 

We would like to show $d(z-p_1+p_2,p)^2-u-w'(p)\geq0$ for all $p\in P$. 
If $p\notin B(p_2,2\sqrt{5}M)$, then Lemma \ref{bddradius} implies 
\begin{align*}
d(z-p_1+p_2,p)&\geq d(p,p_2)-d(z-p_1+p_2,p_2)=d(p,p_2)-d(z,p_1)\\
&>2\sqrt{5}M-\sqrt{5}M=\sqrt{5}M. 
\end{align*}
It follows that $d(z-p_1+p_2,p)^2-u-w'(p)>5M^2-4M^2>0$. 
When $p$ is in $B(p_2,2\sqrt{5}M)$, 
$p+p_1-p_2$ is in $P\cap B(p_1,2\sqrt{5}M)$ and 
\[
d(z-p_1+p_2,p)^2-u-w'(p)=d(z,p+p_1-p_2)^2-u-w(p+p_1-p_2)\geq0. 
\]
Hence $(z-p_1+p_2,u)$ is empty in $(P,w')$, 
and so $\tau'-p_1+p_2$ is in $\Del(P,w')$. 
\end{proof}

Lemma \ref{Iptauw} and Proposition \ref{sliver} are 
the interpretation in our context 
of Claim 14 and the Sliver Theorem of \cite{CDEFT}. 
The following lemma claims that 
if the weight on $p$ is not in $I(p,\tau,\lambda)$, 
then the face $\tau\cup\{p\}$ does not appear in $\Del(P,w)$ and 
that the length of $I(p,\tau,\lambda)$ is small 
when $\tau\cup\{p\}$ is a `thin' simplex. 

\begin{lem}\label{Iptauw}
For any $p\in P$, $\tau\in D_k$ with $\tau\cup\{p\}\in D_{k+1}$ and 
a map $\lambda$ from $\tau$ to $[0,(M/3)^2]$, 
there exists an open interval $I(p,\tau,\lambda)\subset\R$ such that 
the following are satisfied. 
\begin{enumerate}
\item If there exist $w\in\mathcal{W}$ and $\tilde\tau\in\Del(P,w)$ 
such that $\lambda=w|\tau$ and $\tau\cup\{p\}\subset\tilde\tau$, 
then $w(p)$ is in $I(p,\tau,\lambda)$. 
\item The length of $I(p,\tau,\lambda)$ is $8Md(p,H(\tau))$. 
\end{enumerate}
Moreover, the following also holds. 
\begin{enumerate}\setcounter{enumi}{2}
\item Let $p_i\in P$, 
$\tau_i\in D_k$ with $\tau_i\cup\{p_i\}\in D_{k+1}$ and 
$\lambda_i:\tau_i\to[0,(M/3)^2]$ for $i=1,2$. 
If there exists $t\in\R^d$ such that 
$\tau_1+t=\tau_2$, $p_1+t=p_2$ and $\lambda_1(q)=\lambda_2(q+t)$ 
for all $q\in\tau_1$, 
then $I(p_1,\tau_1,\lambda_1)$ is equal to $I(p_2,\tau_2,\lambda_2)$. 
\end{enumerate}
\end{lem}
\begin{proof}
Let $p\in P$, $\tau\in D_k$ with $\tau\cup\{p\}\in D_{k+1}$ and 
$\lambda:\tau\to[0,(M/3)^2]$. 
Put $\tau'=\tau\cup\{p\}$. 
Let $\pi$ be the orthogonal projection from $H(\tau')$ to $H(\tau)$. 
Note that $d(p,H(\tau))$ equals $d(p,\pi(p))$. 

Let $(z,u)$ be the orthosphere of 
the weighted points $(q,\lambda(q))$ for $q\in\tau$. 
Put 
\[
a=d(p,H(\tau))^2+d(z,\pi(p))^2-u\in\R. 
\]
We define 
\[
I(p,\tau,\lambda)=\{x\in\R\mid\lvert x-a\rvert<4Md(p,H(\tau))\}. 
\]
Then it is easy to see that (2) and (3) hold. 

Let us show (1). 
Let $w\in\mathcal{W}$ with $\lambda=w|\tau$. 
Suppose that there exists $\tilde\tau\in\Del(P,w)$ 
which contains $\tau'=\tau\cup\{p\}$. 
Let $(z',u')$ be the orthosphere of $\tau'$ with respect to $w$. 
From Lemma \ref{tau'} and \ref{bddradius}, 
we have $u'<4M^2$. 

By Lemma \ref{tau'}, $\pi(z')=z$. 
Let $b$ be the signed distance from $z$ to $z'$: 
$b$ is positive if and only if 
$z'$ and $p$ lie on the same side of $H(\tau)$ in $H(\tau')$, 
and $b$ is negative if and only if 
$z'$ and $p$ lie on different sides of $H(\tau)$ in $H(\tau')$. 
Then we have 
\begin{align*}
d(z',p)^2
&=d(z,\pi(p))^2+(d(p,\pi(p))-b)^2 \\
&=d(z,\pi(p))^2+(d(p,H(\tau))-b)^2. 
\end{align*}
By Lemma \ref{tau'}, we also have $d(z',z)^2=u'-u$. 
Combining these equations with $d(z',p)^2=u'+w(p)$, 
we get 
\begin{align*}
w(p)&=d(z,\pi(p))^2+(d(p,H(\tau))-b)^2-d(z',z)^2-u \\
&=a-2bd(p,H(\tau)). 
\end{align*}
Moreover, $\lvert b\rvert=d(z',z)=\sqrt{u'-u}\leq\sqrt{u'}<2M$, 
and so $w(p)$ belongs to $I(p,\tau,\lambda)$. 
\end{proof}

The following proposition, 
which will be used in the proof of Proposition \ref{pretriangle}, 
enables us to get rid of slivers. 
Indeed, condition (2) means that 
if the weight on $p$ is not in $I(p,w)$, then 
the resultant weighted Delaunay triangulation does not contain 
the face $\tau\cup\{p\}$ with $d(p,H(\tau))$ small. 

\begin{prop}\label{sliver}
There exists a constant $c_1>0$ depending only on $d$ 
such that the following holds. 

For any $p\in P$ and $w\in{\cal W}$, 
there exists an open subset $I(p,w)\subset\R$ such that 
the following are satisfied. 
\begin{enumerate}
\item $[0,(M/3)^2]\setminus I(p,w)$ is not empty. 
\item For $1\leq k<d$, 
suppose that $\tau\in D_k$ with $\tau\cup\{p\}\in D_{k+1}$ and 
$d(p,H(\tau))\leq c_1M$. 
If $w'\in{\cal W}$ satisfies 
$w'(q)=w(q)$ for all $q\in\tau$ and $w'(p)\notin I(p,w)$, 
then there does not exist $\tilde\tau\in\Del(P,w')$ 
such that $\tau\cup\{p\}\subset\tilde\tau$. 
\end{enumerate}
Moreover, $I(p,w)$ is locally derived in the following sense. 
\begin{enumerate}\setcounter{enumi}{2}
\item If $p_1,p_2\in P$ and $w_1,w_2\in{\cal W}$ satisfy 
\[
(P-p_1)\cap B(0,2\sqrt{5}M)=(P-p_2)\cap B(0,2\sqrt{5}M)
\]
and $w_1(q)=w_2(q-p_1+p_2)$ 
for all $q\in(P\setminus\{p_1\})\cap B(p_1,2\sqrt{5}M)$, 
then $I(p_1,w_1)=I(p_2,w_2)$. 
\end{enumerate}
\end{prop}
\begin{proof}
We establish the claim for the constant 
$c_1=(72\cdot10^{d^2}d)^{-1}$. 

For $p\in P$ and $w\in{\cal W}$, we define 
\[
I(p,w)=\bigcup_{k=1}^{d-1}\bigcup_{\tau}I(p,\tau,w|\tau), 
\]
where the second union runs over all $\tau\in D_k$ such that 
$\tau\cup\{p\}\in D_{k+1}$ and $d(p,H(\tau))\leq c_1M$. 

Let us check (1). 
For any $\tau\in D_k$ such that $\tau\cup\{p\}\in D_{k+1}$, 
we have $\tau\subset P\cap B(p,2\sqrt{5}M)$. 
By Lemma \ref{sepasyn}, 
\[
\#(P\cap B(p,2\sqrt{5}M))\leq(4\sqrt{5}+1)^d<10^d. 
\]
It follows that 
\[
\#\{\tau\in D_k\mid\tau\cup\{p\}\in D_{k+1}\}<10^{(k+1)d}
\]
for each $k=1,2,\dots,d-1$. 
Hence 
\[
\#\left(
\bigcup_{k=1}^{d-1}\{\tau\in D_k\mid\tau\cup\{p\}\in D_{k+1}\}
\right)
<10^{2d}+10^{3d}+\dots+10^{d^2}<d10^{d^2}. 
\]
If $\lvert\cdot\rvert$ denotes the Lebesgue measure on $\R$ 
and if $d(p,H(\tau))\leq c_1M$, 
then by Lemma \ref{Iptauw} (2) 
\[
\lvert I(p,\tau,w|\tau)\rvert\leq8c_1M^2. 
\]
Hence 
\[
\lvert I(p,w)\rvert<8d10^{d^2}c_1M^2=M^2/9, 
\]
and therefore $[0,(M/3)^2]\setminus I(p,w)$ is not empty. 

We now verify (2). 
Suppose that $\tau$ belongs to $D_k$, 
$\tau\cup\{p\}$ belongs to $D_{k+1}$ and 
$w'\in{\cal W}$ satisfies 
$w'(q)=w(q)$ for all $q\in\tau$ and $w'(p)\notin I(p,w)$. 
If $d(p,H(\tau))\leq c_1M$, 
then, by the definition of $I(p,w)$, 
it contains $I(p,\tau,w|\tau)$, 
which is equal to $I(p,\tau,w'|\tau)$. 
It follows that 
$w'(p)$ does not belong to $I(p,\tau,w'|\tau)$. 
Hence, by Lemma \ref{Iptauw} (1), 
there does not exist $\tilde{\tau}\in\Del(P,w')$ 
which contains $\tau\cup\{p\}$. 

Finally, let us consider (3). 
Suppose that $p_1,p_2\in P$ and $w_1,w_2\in{\cal W}$ satisfy 
\[
(P-p_1)\cap B(0,2\sqrt{5}M)=(P-p_2)\cap B(0,2\sqrt{5}M)
\]
and $w_1(q)=w_2(q-p_1+p_2)$ 
for all $q\in(P\setminus\{p_1\})\cap B(p_1,2\sqrt{5}M)$. 
By symmetry, 
it suffices to show $I(p_1,w_1)\subset I(p_2,w_2)$. 
If $x\in I(p_1,w_1)$, then by definition 
there exists $\tau\in D_k$ such that 
$\tau\cup\{p_1\}\in D_{k+1}$, $d(p_1,H(\tau))\leq c_1M$ 
and $x\in I(p_1,\tau,w_1|\tau)$. 
By Lemma \ref{Dkloc}, 
$\tau'=\tau-p_1+p_2$ belongs to $D_k$ and 
$\tau'\cup\{p_2\}$ belongs to $D_{k+1}$. 
It follows from Lemma \ref{Iptauw} (3) that 
$I(p_1,\tau,w_1|\tau)$ is equal to $I(p_2,\tau',w_2|\tau')$. 
We still have $d(p_2,H(\tau'))\leq c_1M$, 
because $\tau'$ is just a translation of $\tau$. 
Therefore $I(p_2,\tau',w_2|\tau')$ is contained in $I(p_2,w_2)$, 
which implies $x\in I(p_2,w_2)$. 
This completes the proof. 
\end{proof}

\section{$\phi$-regular triangulations of $\R^d$}

In this section, 
we would like to apply the argument of the last section 
to countable subsets of $\R^d$ 
arising from a minimal dynamical system on a Cantor set. 
First, we would like to recall 
the notion of $\phi$-regularity and local derivedness 
(\cite[Section 4]{GMPS1}). 
Let $X\subset\Omega$ be a flat Cantor transversal 
of a free minimal action $\phi$ of $\R^d$ on $\Omega$ 
and let $R_\phi\subset X\times X$ be the \'etale equivalence relation 
induced from $(\Omega,\phi)$. 

Suppose that for each $x$ in $X$, we have a subset $P(x)$ of $\R^d$. 
We say that this collection is $\phi$-regular if the following hold. 
\begin{enumerate}
\item For any $x,\phi^p(x)\in X$, $P(\phi^p(x)) = P(x)-p$. 
\item If $x$ is in $X$ and $K\subset\R^d$ is compact, 
then there exists a neighbourhood $U$ of $x$ in $X$ such that 
\[
P(x')\cap K=P(x)\cap K
\]
for all $x'$ in $U$. 
\end{enumerate}

For a clopen subset $U\subset X$, 
we define 
\[
P_U(x)=\{p\in\R^d\mid\phi^p(x)\in U\}. 
\]
In what follows, to simplify notation, 
we will often denote the family $\{P_U(x)\}_{x\in X}$ by $P_U$. 
We remark that, 
from the definition of flat Cantor transversals, 
there exist positive real numbers $M_0,M_1>0$ such that 
$P_X(x)$ is $M_0$-separated and $M_1$-syndetic for every $x\in X$. 
The following result is 
an easy consequence of the definition and we omit the proof. 

\begin{lem}
Let $U$ be a clopen subset of $X$. 
The family of sets $P_U(x)$ for $x\in X$ is $\phi$-regular. 
Conversely, if $\{P(x)\}_{x\in X}$ is a $\phi$-regular family, then 
\[
U=\{x\in X\mid0\in P(x)\}
\]
is clopen. 
\end{lem}

We consider a family $\{\mathcal{T}(x)\}_{x\in X}$ 
of tessellations of $\R^d$ 
which are indexed by the points of $X$. 
We say that this collection is $\phi$-regular 
if the following hold. 
\begin{enumerate}
\item For any $x,\phi^p(x)\in X$, 
$\mathcal{T}(\phi^p(x))=\mathcal{T}(x)-p$. 
\item If $x$ is in $X$ and $t$ is in $\mathcal{T}(x)$, 
then there is a neighbourhood $U$ of $x$ such that 
$t$ is in $\mathcal{T}(x')$ for all $x'$ in $U$. 
\end{enumerate}

If $P$ is a $\varphi$-regular family, 
we say that it is $M$-syndetic, ($M$-separated, respectively) 
if $P(x)$ is $M$-syndetic ($M$-separated, respectively) 
for each $x$ in $X$. 

Let $P(x),P'(x)$, $x\in X$ be two families of subsets of $\R^d$. 
We say that $P'$ is locally derived from $P$ 
if there is a constant $R>0$ such that, 
for any $x_1,x_2\in X$ and $u_1,u_2\in\R^d$, 
if $u_1$ is in $P'(x_1)$ and 
\[
(P(x_1)-u_1)\cap B(0,R)=(P(x_2)-u_2)\cap B(0,R), 
\]
then $u_2$ is in $P'(x_2)$. 
In a similar way, we extend this definition 
replacing either $P$, $P'$ or both with families of tessellations. 
For example, a family of tessellations $\{\mathcal{T}(x)\}_{x\in X}$ is 
said to be locally derived from another family of tessellations 
$\{\mathcal{T}'(x)\}_{x\in X}$,  
if there exists $R>0$ such that 
for any $x_1,x_2\in X$ and $u_1,u_2\in\R^d$, 
if $t\in\mathcal{T}'(x_1)-u_1$ contains the origin and 
\[
\{s\in\mathcal{T}(x_1)-u_1\mid s\cap B(0,R)\neq\emptyset\}
=\{s\in\mathcal{T}(x_2)-u_2\mid s\cap B(0,R)\neq\emptyset\}, 
\]
then $t$ is in $\mathcal{T}'(x_2)-u_2$. 
The following result is easily derived 
from the definitions and we omit the proof. 

\begin{lem}\label{localderive}
If $P$ is a $\phi$-regular family and $P'$ is locally derived from $P$, then 
$P'$ is also $\phi$-regular. 
Analogous statements hold 
replacing $P$, $P'$ or both with families of tessellations. 
\end{lem}

\bigskip

Next, we turn to the issue of the existence of 
$\phi$-regular, separated and syndetic sets for flat Cantor transversals. 

\begin{lem}\label{aperiodic}
For any $M>0$, 
there exists a clopen partition $\{V_1,V_2,\dots,V_n\}$ of $X$ 
such that $V_i\cap\phi^p(V_i)=\emptyset$ 
for all nonzero $p\in B(0,M)$ and $i=1,2,\dots,n$. 
\end{lem}
\begin{proof}
For each $x\in X$, let $E_x=\{p\in B(0,M)\mid\phi^p(x)\in X\}$. 
By Definition \ref{Ctrans} (4), 
there exists a clopen neighbourhood $U_x\subset X$ of $x$ such that 
$E_x=E_y$ for all $y\in U_x$. 
By Definition \ref{Ctrans} (3), $E_x$ is a finite set, 
and so we may assume that 
$U_x$ is chosen sufficiently small so that 
$U_x\cap\phi^p(U_x)=\emptyset$ for every nonzero $p\in E_x$. 
Therefore we have 
$U_x\cap\phi^p(U_x)=\emptyset$ for all nonzero $p\in B(0,M)$. 

By the compactness of $X$, 
we can select a finite set $x_1,x_2,\dots,x_n\in X$ 
such that $U_{x_i}$'s cover $X$. 
Put 
\[
V_k=U_{x_k}\setminus\bigcup_{i=1}^{k-1}U_{x_i}. 
\]
Then $\{V_1,V_2,\dots,V_n\}$ is a clopen partition of $X$ and 
$V_k\cap\phi^p(V_k)=\emptyset$ for all nonzero $p\in B(0,M)$. 
\end{proof}

By using the lemma above, we can prove the following 
in a similar way to \cite[Proposition 4.4]{GMPS1}. 
We omit the proof. 

\begin{lem}\label{M2M}
Suppose that $P_X$ is $M_1$-syndetic. 
For any $M>M_1$, 
there exists a clopen subset $U\subset X$ 
such that $P_U$ is $M$-separated and $2M$-syndetic. 
\end{lem}

Now we are ready to prove Proposition \ref{pretriangle}, 
which claims the existence of a weight function $w$ such that 
any $d$-simplex $\tau$ in the resultant $\Del(\cdot,\cdot)$ 
is not so `thin' (i.e. for any $\tau\in\Del(P,w)$ and $p\in\tau$, 
the distance from $p$ to the $(d-1)$-dimensional affine subspace 
containing $\tau\setminus\{p\}$ is bounded below). 
We define a positive constant $c_2$ by 
\[
c_2=\frac{c_1^{d-1}}{(d-1)!V_{d-1}(\sqrt{5})^{d-1}}, 
\]
where $c_1$ is the constant obtained in Proposition \ref{sliver} 
and $V_{d-1}$ is the volume of $(d-1)$-dimensional unit ball. 
We remark that $c_2$ depends only on $d$. 

\begin{prop}\label{pretriangle}
Suppose that $P_X$ is $M_1$-syndetic. 
For any $M>M_1$, 
there exist a clopen subset $U\subset X$ and 
a locally constant function $w:U\to[0,(M/3)^2]$ 
such that the following hold, 
where $w_x:P_U(x)\to[0,(M/3)^2]$ is defined by $w_x(p)=w(\phi^p(x))$. 
\begin{enumerate}
\item $P_U$ is $M$-separated and $2M$-syndetic. 
\item For any $x\in X$ and $\tau\in\Del(P_U(x),w_x)$, 
$\tau$ is contained in a closed ball with radius less than $\sqrt{5}M$. 
\item For any $x\in X$ and $\tau\in\Del(P_U(x),w_x)$, 
the volume of $\Conv(\tau)$ is greater than $c_1^{d-1}M^d/d!$. 
\item For any $x\in X$, $\tau\in\Del(P_U(x),w_x)$ and $p\in\tau$, 
$d(p,H(\tau\setminus\{p\}))$ is greater than $c_2M$. 
\end{enumerate}
\end{prop}
\begin{proof}
By Lemma \ref{M2M}, 
there exists a clopen subset $U\subset X$ 
such that $P_U$ is $M$-separated and $2M$-syndetic. 
We would like to construct $w:U\to[0,(M/3)^2]$ inductively. 
By Lemma \ref{aperiodic}, 
we can find a clopen partition $\{V_1,V_2,\dots,V_n\}$ of $X$ 
such that $V_i\cap\phi^p(V_i)=\emptyset$ 
for all nonzero $p\in B(0,2\sqrt{5}M)$ and $i=1,2,\dots,n$. 
Let $w_0$ be the constant function on $U$ with value zero. 

Suppose that 
we have fixed a locally constant function $w_{i-1}:U\to[0,(M/3)^2]$. 
Take $x\in X$ and 
consider the $M$-separated and $2M$-syndetic subset $P_U(x)$ 
with the weight function $w_{i-1,x}(p)=w_{i-1}(\phi^p(x))$.  
Let us define the map $w_{i,x}:P_U(x)\to[0,(M/3)^2]$. 
For $p\in P_U(x)$ with $\phi^p(x)\notin U\cap V_i$, 
we put $w_{i,x}(p)=w_{i-1,x}(p)$. 
Suppose that $p\in P_U(x)$ satisfies $\phi^p(x)\in U\cap V_i$. 
We define $w_{i,x}(p)$ to be the minimum value of 
the closed set $[0,(M/3)^2]\setminus I(p,w_{i-1,x})$, 
which is not empty from Proposition \ref{sliver} (1). 
Let $\omega_i$ denote the function given for $x\in U$, 
by $\omega_i(x)=\omega_{i,x}(0)$. 
Then $\omega_i=\omega_{i-1}$ on $U\setminus V_i$ and 
therefore $\omega_i$ is locally constant by assumption. 
If $x,y\in U\cap V_i$ are close enough, then 
we have $P_U(x)\cap B(0,2\sqrt{5}M)=P_U(y)\cap B(0,2\sqrt{5}M)$. 
Moreover, for any non-zero $q\in P_U(x)\cap B(0,2\sqrt{5}M)$, 
by the choice of $V_i$, $\phi^q(x)$ is not in $V_i$. 
Likewise $\phi^q(y)$ is not in $V_i$. 
Hence $w_{i,x}(q)=w_{i-1,x}(q)$ equals $w_{i,y}(q)=w_{i-1,y}(q)$, 
if $x$ and $y$ are sufficiently close. 
By Proposition \ref{sliver} (3), for such $x,y$, 
we can conclude that $I(0,w_{i-1,x})$ equals $I(0,w_{i-1,y})$, 
which implies $w_{i,x}(0)=w_{i,y}(0)$. 
Thus, the function $w_i$ is locally constant. 
Repeating this procedure, we get $w_n$. 
Put $w=w_n$. 

Take $x\in X$ and let $\tau=\{p_0,p_1,\dots,p_d\}\in\Del(P_U(x),w_x)$. 
Let us show (2). 
Let $z$ be the orthocentre of $\tau$ with respect to $w_x$. 
By Lemma \ref{bddradius}, $d(z,p_i)$ is less than $\sqrt{5}M$. 
It follows that $\tau$ is contained in $B(z,\sqrt{5}M)$. 

We next verify (3). 
By Lemma \ref{bddradius}, 
for any $p_i,p_j\in\tau$, $d(p_i,p_j)$ is less than $2\sqrt{5}M$. 
Hence, from the choice of clopen sets $V_l$'s, 
$\phi^{p_i}(x)$ and $\phi^{p_j}(x)$ are not contained 
in the same $V_l$ for $i\neq j$. 
Without loss of generality, 
we may assume $\phi^{p_i}(x)\in V_{l_i}$ and $l_0<l_1<\dots<l_d$. 
Let $\tau_k=\{p_0,p_1,\dots,p_k\}$ for each $k=1,2,\dots,d-1$. 
Since $P_U$ is $M$-separated, 
the one-dimensional volume, that is, 
the length of $\Conv(\tau_1)$ is not less than $M$. 

For each $k=2,3,\dots,d$, 
we would like to see $d(p_k,H(\tau_{k-1}))>c_1M$. 
Suppose that $d(p_k,H(\tau_{k-1}))$ is not greater than $c_1M$. 
By definition of $w$, 
for all $q\in\tau_{k-1}$, 
$w_x(q)$ is equal to $w_{l_k-1,x}(q)$. 
In addition, 
$w_x(p_k)=w_{l_k,x}(p_k)$ does not belong to $I(p_k,w_{l_k-1,x})$. 
By applying Proposition \ref{sliver} (2) to 
$p_k$, $w_{l_k-1,x}$, $\tau_{k-1}$ and $w_x$, 
we can conclude that 
there does not exist $\tilde{\tau}\in\Del(P_U(x),w_x)$ 
such that $\tau_{k-1}\cup\{p_k\}\subset\tilde{\tau}$. 
This is clearly a contradiction, and so 
$d(p_k,H(\tau_{k-1}))>c_1M$ for all $k=2,3,\dots,d$. 
It follows that 
the $d$-dimensional volume of $\Conv(\tau_d)$ is greater than 
\[
M\times\frac{1}{2\cdot3\cdot\dots\cdot d}(c_1M)^{d-1}
=\frac{c_1^{d-1}M^d}{d!}. 
\]

We next verify (4). 
By the argument above, for any $p\in\tau$, 
$\tau\setminus\{p\}$ is contained 
in a $(d-1)$-dimensional closed disk with radius less than $\sqrt{5}M$. 
It follows that 
$(d-1)$-dimensional volume of $\Conv(\tau\setminus\{p\})$ is 
less than $V_{d-1}(\sqrt{5}M)^{d-1}$, 
where $V_{d-1}$ is the volume of $(d-1)$-dimensional unit ball. 
Hence 
\[
d(p,H(\tau\setminus\{p\}))
>\frac{dc_1^{d-1}M^d}{d!V_{d-1}(\sqrt{5}M)^{d-1}}
=c_2M
\]
for all $p\in\tau$. 
\end{proof}

By the proposition above, 
we have obtained $\Del(P_U(x),w_x)$ containing no slivers. 
However, some extra work is necessary to obtain triangulations of $\R^d$, 
because $\{\Conv(\tau)\mid\tau\in\Del(P_U(x),w_x)\}$ may not be 
a triangulation in general, or equivalently 
$\#\delta(z,u)$ may not equal $d{+}1$ 
for $(z,u)\in\mathcal{Z}(P_U(x),w_x)$. 
We would like to triangulate $\Conv(\delta(z,u))$ 
(for $(z,u)$ with $\#\delta(z,u)>d{+}1$) 
in a locally derived manner. 
We need an elementary fact about triangulations of convex polytopes. 
In what follows, we use the notation in \cite[Lemma 1.1]{St}. 
But the reader should be warned that 
Lemma 1.1 of \cite{St} contains only the statement. 
In page 160 of \cite{Mat}, 
one can find a construction called bottom-vertex triangulation, 
which is the same as that used here. 

Equip $\R^d$ with the lexicographic ordering. 
Namely, for $p=(p_1,p_2,\dots,p_d)$ and $q=(q_1,q_2,\dots,q_d)$ 
in $\R^d$, $p$ is less than $q$, if there exists $i$ such that 
$p_i<q_i$ and $p_j=q_j$ for all $j<i$. 
Note that the lexicographic ordering is invariant under translation 
in the sense that $p<q$ implies $p+r<q+r$ for any $r\in\R^d$. 

Let $K\subset\R^d$ be a convex polytope. 
For each non-empty face $F$ of $K$, 
we define $m(F)$ to be the minimum element 
in the set of vertices (i.e. extremal points) of $F$. 
A $(d+1)$-tuple $\Phi=(F_0,F_1,\dots,F_d)$ is called a flag of faces, 
if $F_i$ is an $i$-dimensional face of $K$ 
and $F_0\subset F_1\subset\dots\subset F_d=K$. 
Call $\Phi$ a full flag, 
if $m(F_i)$ is not in $F_{i-1}$ for $1\leq i\leq d$. 
For a full flag $\Phi$, 
we put $\tau_\Phi=\{F_0,m(F_1),\dots,m(F_d)\}$. 
Then 
\[
\{\Conv(\tau_\Phi)\mid\Phi\text{ is a full flag of }K\}
\]
is a triangulation of $K$, 
in which the cells meet face to face. 
Thus, the union of those $\Conv(\tau_\Phi)$'s is equal to $K$ 
and for any full flags $\Phi$ and $\Phi'$, 
we have 
\[
\Conv(\tau_\Phi)\cap\Conv(\tau_{\Phi'})
=\Conv(\tau_\Phi\cap\tau_{\Phi'}). 
\]

We are able to apply this argument to tessellations of $\R^d$ 
by convex polytopes. 
Let $\mathcal{T}$ be a tessellation of $\R^d$ 
by convex polytopes. 
Suppose that cells in $\mathcal{T}$ meet face to face. 
Let $V$ be the set of all vertices of all cells in $\mathcal{T}$. 
For each cell in $\mathcal{T}$, 
by considering all full flags of it, 
we obtain its triangulation. 
Let $\mathcal{D}$ be the union of them, i.e. 
\[
\mathcal{D}=\{\Conv(\tau_\Phi)\mid
\Phi\text{ is a full flag of some cell }t\in\mathcal{T}\}. 
\]
Notice that $\tau_\Phi$ is in $\Delta_d(V)$. 
It is not so hard to see the following. 
\begin{itemize}
\item Each $\Conv(\tau)$ in $\mathcal{D}$ is contained 
in a cell of $\mathcal{T}$. 
\item $\mathcal{D}$ is a triangulation of $\R^d$. 
\item For any $\Conv(\tau),\Conv(\tau')\in\mathcal{D}$, 
\[
\Conv(\tau)\cap\Conv(\tau')=\Conv(\tau\cap\tau'), 
\]
that is, cells of $\mathcal{D}$ meet face to face. 
\end{itemize}
Moreover, we can make the following remark. 
Since the lexicographic ordering on $\R^d$ is invariant by translation, 
$\mathcal{D}$ is locally derived from $\mathcal{T}$. 
Thus, if $\mathcal{T}$ is a $\phi$-regular family of tessellations, 
then $\mathcal{D}$ also becomes $\phi$-regular 
by Lemma \ref{localderive}. 

From Lemma \ref{Laguerre}, Proposition \ref{pretriangle} and 
the discussion above, 
we then prove the following proposition. 

\begin{prop}\label{triangle}
Suppose that $P_X$ is $M_1$-syndetic. 
For any $M>M_1$, 
there exist a clopen subset $U\subset X$ and 
a family of subsets $\mathcal{D}(x)\subset\Delta_d(P_U(x))$ 
for $x\in X$ satisfying each of the following conditions 
for all $x\in X$. 
\begin{enumerate}
\item $P_U(x)$ is $M$-separated and $2M$-syndetic. 
\item The collection $\{\Conv(\tau)\mid\tau\in\mathcal{D}(x)\}$ is 
a tessellation of $\R^d$. 
\item For any $\tau,\tau'\in\mathcal{D}(x)$, 
we have $\Conv(\tau)\cap\Conv(\tau')=\Conv(\tau\cap\tau')$. 
\item For any $p\in P_U(x)$, 
there exists $\tau\in\mathcal{D}(x)$ such that $p\in\tau$. 
\item Any $\tau\in\mathcal{D}(x)$ is contained 
in a closed ball with radius less than $\sqrt{5}M$. 
\item For any $\tau\in\mathcal{D}(x)$ and $p\in\tau$, 
$d(p,H(\tau\setminus\{p\}))$ is greater than $c_2M$. 
\item For any $p\in P_X(x)$, 
$\tau$ is in $\mathcal{D}(x)$ if and only if 
$\tau-p$ is in $\mathcal{D}(\phi^p(x))$. 
\item For any $\tau\in\mathcal{D}(x)$, 
there exists an open neighbourhood $V\subset X$ of $x\in X$ 
such that $\tau\in\mathcal{D}(y)$ for all $y\in V$. 
\end{enumerate}
\end{prop}
\begin{proof}
From Proposition \ref{pretriangle}, 
we get a clopen subset $U\subset X$ and 
a locally constant function $w:U\to[0,(M/3)^2]$ 
satisfying the conditions given there. 
The first condition follows at once 
from Proposition \ref{pretriangle} (1). 
For each $x\in X$, 
we let $w_x$ denote the map from $P_U(x)$ to $[0,(M/3)^2]$ 
defined by $w_x(p)=w(\phi^p(x))$. 

By Lemma \ref{Laguerre}, for each $x\in X$, 
\[
\{\Conv(\delta(z,u))\mid(z,u)\in\mathcal{Z}(P_U(x),w_x)\}
\]
is a tessellation of $\R^d$ by convex polytopes, 
in which cells meet face to face. 
It follows from the discussion above that 
there exists $\mathcal{D}(x)\subset\Delta_d(P_U(x))$ such that 
\begin{itemize}
\item[(i)] Each $\tau$ in $\mathcal{D}(x)$ is contained 
in some $\delta(z,u)$ in $\mathcal{Z}(P_U(x),w_x)$. 
\item[(ii)] The collection $\{\Conv(\tau)\mid\tau\in\mathcal{D}(x)\}$ is 
a tessellation of $\R^d$. 
\item[(iii)] For any $\tau,\tau'\in\mathcal{D}(x)$, 
we have $\Conv(\tau)\cap\Conv(\tau')=\Conv(\tau\cap\tau')$. 
\end{itemize}
From (ii) and (iii), the conditions (2) and (3) are immediate. 
The condition (4) is clear. 
From (i), we see that $\tau$ belongs to $\Del(P_U(x),w_x)$, 
and so the conditions (5) and (6) follow from 
Proposition \ref{pretriangle}. 

The conditions (7) and (8) easily follow from the continuity of $w$ 
and the local derivedness of the construction of $\mathcal{D}(x)$. 
\end{proof}

\section{Well-separated tessellations of $\R^d$}

In the last section, 
we constructed a triangulation of $\R^d$ with several nice properties. 
In this section, 
we study how to construct a tessellation 
which is `dual' to the triangulation, that is, 
the vertices in the triangulation become cells and 
cells in the triangulation become vertices. 
The tessellation will possess various nice properties. 
Before beginning, 
we make the following remark. 
In the actual application, 
we will begin with a $\phi$-regular collection $P(x)$, 
find weighted Delaunay triangulations without slivers, 
and construct a collection of tessellations $\mathcal{T}(x)$ on $\R^d$. 
It is worth noting as we proceed, that 
all of our construction are `locally derived' 
in the appropriate sense, 
and so the resulting collection $\mathcal{T}(x)$'s is $\phi$-regular 
by application of Lemma \ref{localderive}. 

\begin{df}
Let $\mathcal{T}$ be a tessellation of $\R^d$. 
\begin{enumerate}
\item We say that $\mathcal{T}$ has capacity $C>0$, 
if each element of $\mathcal{T}$ contains an open ball of radius $C$. 
\item We say that $\mathcal{T}$ is $K$-separated for $K>0$, 
if the following hold. 
\begin{enumerate}
\item For any $n\in\N$ and $t_0,t_1,\dots,t_n\in\mathcal{T}$, 
if there exists $p\in\R^d$ such that 
$d(t_i,p)<K$ for all $i=0,1,\dots,n$, 
then $t_0\cap t_1\cap\dots\cap t_n$ is not empty. 
\item Any distinct $d+2$ elements of $\mathcal{T}$ have trivial intersection. 
\end{enumerate}
\item We let the diameter of $\mathcal{T}$ to be 
the supremum of the diameter of its elements 
and denote it by $\diam(\mathcal{T})$. 
When $\{\mathcal{T}(x)\mid x\in X\}$ is a $\phi$-regular family of 
tessellations, then for any $x,y\in X$ and $t\in\mathcal{T}(x)$, 
it follows from the minimality of $\phi$ that 
there exists $p\in P_X(y)$ such that $t\in\mathcal{T}(\phi^p(y))$. 
Therefore $\diam(\mathcal{T}(x))$ does not depend on $x\in X$ 
and we denote this value by $\diam(\mathcal{T})$. 
\end{enumerate}
\end{df}

Suppose that 
we are given an $M$-separated and $2M$-syndetic subset $P$ of $\R^d$ 
and a subset $\mathcal{D}\subset\Delta_d(P)$ satisfying the following. 
\begin{itemize}
\item The collection $\{\Conv(\tau)\mid\tau\in\mathcal{D}\}$ 
is a tessellation of $\R^d$. 
\item For any $\tau,\tau'\in\mathcal{D}$, we have 
$\Conv(\tau)\cap\Conv(\tau')=\Conv(\tau\cap\tau')$. 
\item For any $p\in P$, 
there exists $\tau\in\mathcal{D}$ such that $p\in\tau$. 
\item Any $\tau\in\mathcal{D}$ is contained 
in a closed ball with radius less than $\sqrt{5}M$. 
\item For any $\tau\in\mathcal{D}$ and $p\in\tau$, 
the distance from $p$ to $H(\tau\setminus\{p\})$ is 
greater than $c_2M$. 
\end{itemize}

For $p\in P$, 
we denote $\{\tau\in\mathcal{D}\mid p\in\tau\}$ by $\mathcal{D}_p$. 
For every $p\in P$ and $\tau\in\mathcal{D}_p$, 
we define 
\[
T(p,\tau)=\left\{\sum_{q\in\tau}\lambda_qq\mid
\sum_q\lambda_q=1,\lambda_p\geq\lambda_q\geq0
\text{ for all }q\in\tau\right\}. 
\]
It is easily seen that 
the collection $\{T(p,\tau)\mid p\in P,\tau\in\mathcal{D}_p\}$ 
is a tessellation of $\R^d$. 
For any $p\in P$, we set 
\[
T(p)=\bigcup_{\tau\in\mathcal{D}_p}T(p,\tau). 
\]
Again the collection $\mathcal{T}=\{T(p)\mid p\in P\}$ is 
a tessellation of $\R^d$. 
We would like to show that 
this tessellation $\mathcal{T}$ has various nice properties. 
Notice that $T(p)$ may not be convex, but is a union of convex polytopes. 

\begin{figure}[h]
\begin{center}
\includegraphics[height=0.2\textheight]{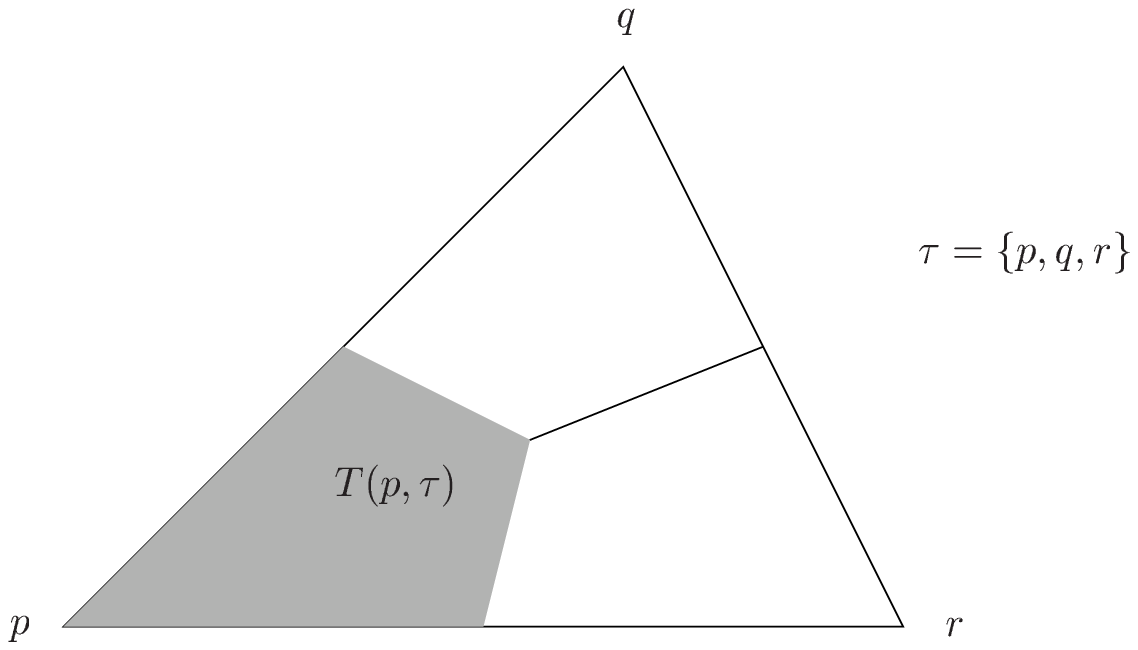}
\caption{$T(p,\tau)$}
\label{fig2}
\end{center}
\end{figure}

For every $p\in P$, we set 
\[
S(p)=\bigcup_{\tau\in\mathcal{D}_p}\Conv(\tau). 
\]
Note that $S(p)$ is homeomorphic to a closed unit ball of $\R^d$ and, 
for any $q\in S(p)$ and $0\leq\lambda\leq1$, 
$\lambda(q-p)+p$ is in $S(p)$. 
It is also easy to see that 
the interior of $S(p)$ is equal to 
\[
\bigcup_{\tau\in\mathcal{D}_p}
\Conv(\tau)\setminus\Conv(\tau\setminus\{p\}). 
\]
Besides, $B(p,c_2M)$ is contained in $S(p)$, 
because the distance from $p$ to $H(\tau\setminus\{p\})$ is 
greater than $c_2M$. 

\begin{lem}\label{capacity}
For every $p\in P$, $B(p,c_2M/(d+1))$ is contained in $T(p)$. 
In particular, $\mathcal{T}$ has capacity $c_2M/(d+1)$. 
\end{lem}
\begin{proof}
Clearly $T(p,\tau)$ contains $(d+1)^{-1}(\Conv(\tau)-p)+p$. 
It follows that $T(p)$ contains $(d+1)^{-1}(S(p)-p)+p$. 
Hence $T(p)$ contains $B(p,c_2M/(d+1))$. 
\end{proof}

Next, we would like to show that 
$\mathcal{T}$ is well-separated. 

\begin{lem}
For every $p\in P$ and $0<\lambda<1$, 
the distance from $\lambda(S(p)-p)+p$ to the boundary of $S(p)$ 
is greater than $(1-\lambda)c_2M$. 
\end{lem}
\begin{proof}
When $d=1,2$, the assertion is an easy geometric observation. 
Let us consider the case of $d>2$. 
Suppose that the minimum distance is achieved 
by $q_1$ in $\lambda(S(p)-p)+p$ and $q_2$ in the boundary of $S(p)$. 
Let $H$ be a two-dimensional plane containing $q_1,q_2$ and $p$. 
For each $\tau\in\mathcal{D}_p$, 
$H\cap\Conv(\tau)$ is either $\{p\}$, 
a line segment with an endpoint $p$ 
or a triangle with a vertex $p$. 
Therefore $H\cap S(p)$ is exactly the same shape 
as in the case of $d=2$. 
It follows that $d(q_1,q_2)$ is greater than $(1-\lambda)c_2M$. 
\end{proof}

For a subset $A\subset\R^d$ and $r>0$, 
we let $B(A,r)$ denote the $r$-neighbourhood of $A$, 
namely, $B(A,r)=\{p\in\R^d\mid d(p,A)<r\}$. 

\begin{lem}\label{neighbour}
For every $p\in P$, 
$B(T(p),c_2M/(d+1))$ is contained in the interior of $S(p)$. 
\end{lem}
\begin{proof}
For each $\tau\in\mathcal{D}_p$, 
we consider  $d(d+1)^{-1}(\Conv(\tau)-p)+p$, which is equal to 
\[
\left\{\sum_{q\in\tau}\lambda_qq\mid
\sum_{q\in\tau}\lambda_q=1,\lambda_p\geq(d+1)^{-1},0\leq\lambda_q\leq1
\text{ for all }q\in\tau\right\}. 
\]
Clearly this set contains $T(p,\tau)$. 
It follows that $T(p)$ is contained in $d(d+1)^{-1}(S(p)-p)+p$. 
By the lemma above, 
the $(d+1)^{-1}c_2M$-neighbourhood of $d(d+1)^{-1}(S(p)-p)+p$ is 
contained in the interior of $S(p)$. 
Therefore we can conclude that 
$B(T(p),(d+1)^{-1}c_2M)$ is contained in the interior of $S(p)$. 
\end{proof}

\begin{lem}\label{separated}
The tessellation $\mathcal{T}$ is $c_2M/(d+1)$-separated. 
\end{lem}
\begin{proof}
Suppose that we are given 
$n\in\N$, $p_0,p_1,\dots,p_n\in P$ and $p\in\R^d$ 
such that $d(T(p_i),p)<c_2M/(d+1)$ for all $i=0,1,\dots,n$. 
From Lemma \ref{neighbour}, 
$p$ is in the interior of $S(p_i)$. 
It follows that there exists $\tau_i\in\mathcal{D}_{p_i}$ such that 
$p$ is in $\Conv(\tau_i)\setminus\Conv(\tau_i\setminus\{p_i\})$. 

For each $i=0,1,\dots,n$, 
$p$ is written as 
\[
p=\sum_{q\in\tau_i}\lambda_{i,q}q, 
\]
where $0\leq\lambda_{i,q}\leq1$ and $\sum_{q\in\tau_i}\lambda_{i,q}=1$, 
and this expression is unique. 
For each $i=1,2,\dots,n$, $p$ belongs to $\Conv(\tau_0)\cap\Conv(\tau_i)$, 
which is equal to $\Conv(\tau_0\cap\tau_i)$. 
If $\tau_0$ does not contain $p_i$, then 
$\tau_0\cap\tau_i$ does not, neither. 
Hence $p$ should belong to $\Conv(\tau_i\setminus\{p_i\})$, 
which is a contradiction. 
Thus, we can conclude that 
$\tau_0$ contains $p_i$ for all $i=0,1,\dots,n$. 
By symmetry, each $\tau_i$ contains all the points $\{p_0,p_1,\dots,p_n\}$. 
Therefore, each $T(p_i,\tau_i)$ contains 
the barycentre $z=(n+1)^{-1}\sum_{i=0}^np_i$. 
It follows that each $T(p_i)$ contains $z$. 

In addition, since the cardinality of $\tau_0$ is $d+1$, 
$n$ must be less than $d+1$. 
It follows that 
any $d+2$ distinct elements have trivial intersection. 
\end{proof}

The following lemma will be necessary to show that 
`boundary points' have measure zero in a later section. 
For $A\subset\R^d$, we let $\partial A$ denote the boundary of $A$. 

\begin{lem}\label{boundary}
For every $p\in P$ and $L>0$, 
the $d$-dimensional volume of $B(\partial T(p),L)$ is 
less than $10^{d^2}2LV_{d-1}(\sqrt{5}M+L)^{d-1}$, 
where $V_{d-1}$ is the $(d-1)$-dimensional volume of 
the unit ball of $\R^{d-1}$. 
\end{lem}
\begin{proof}
For each $\tau\in\mathcal{D}_p$ and 
a bijection $f:\{1,2,\dots,d\}\to\tau\setminus\{p\}$, 
we define 
\[
\kappa(\tau,f)=\left\{
\frac{1}{k+1}\left(p+\sum_{i=1}^kf(i)\right)\mid
k=1,2,\dots,d\right\}. 
\]
Then $\kappa(\tau,f)$ consists of 
$d$ distinct points in $\Conv(\tau)$. 
In particular, $\kappa(\tau,f)$ is contained 
in a $(d-1)$-dimensional disk of radius $\sqrt{5}M$, 
because $\tau$ is contained in a closed ball 
with radius less than $\sqrt{5}M$. 

It is easy to see 
\[
T(p,\tau)=\bigcup\Conv(\{p\}\cup\kappa(\tau,f)), 
\]
where the union runs over all 
bijections $f:\{1,2,\dots,d\}\to\tau\setminus\{p\}$. 
Moreover, we have 
\[
\partial T(p)\cap T(p,\tau)=\bigcup\Conv(\kappa(\tau,f)), 
\]
and so 
\[
B(\partial T(p),L)=\bigcup B(\Conv(\kappa(\tau,f)),L), 
\]
where the union runs over all the pairs of 
$\tau\in\mathcal{D}_p$ and 
bijections $f:\{1,2,\dots,d\}\to\tau\setminus\{p\}$. 
Since $\tau$ is contained in a closed ball 
with radius less than $\sqrt{5}M$, 
for any $q\in\tau$, one has $d(p,q)<2\sqrt{5}M$. 
Since $P$ is $M$-separated, by Lemma \ref{sepasyn} (1), 
we have $\#(P\cap B(p,2\sqrt{5}M))\leq(4\sqrt{5}+1)^d<10^d$. 
The map sending $(\tau,f)$ to $(f(1),f(2),\dots,f(d))$ is injective, 
and so the number of pairs $(\tau,f)$ is less than 
$\#(P\cap B(p,2\sqrt{5}M))^d<(10^d)^d=10^{d^2}$. 

It remains for us to estimate 
the $d$-dimensional volume of $B(\Conv(\kappa(\tau,f)),L)$. 
As mentioned above, $\kappa(\tau,f)$ is contained 
in a $(d-1)$-dimensional disk of radius $\sqrt{5}M$. 
Hence the $d$-dimensional volume of 
$B(\Conv(\kappa(\tau,f)),L)$ is not greater than 
\[
2LV_{d-1}(\sqrt{5}M+L)^{d-1}, 
\]
where $V_{d-1}$ is the $(d-1)$-dimensional volume of 
the unit ball of $\R^{d-1}$. 
This completes the proof. 
\end{proof}

\bigskip

Now we would like to apply the construction of $\mathcal{T}$ 
to the family $\mathcal{D}(x)$ obtained in the last section. 
Let $X\subset\Omega$ be a flat Cantor transversal 
of a free minimal action $\phi$ of $\R^d$ on $\Omega$ 
and let $R_\phi\subset X\times X$ be the \'etale equivalence relation 
induced from $(\Omega,\phi)$. 

\begin{prop}\label{wellsepa}
Suppose that $P_X$ is $M_1$-syndetic. 
For any $M>M_1$, 
there exist a clopen subset $U\subset X$, 
a $\phi$-regular family of tessellations $\mathcal{T}(x)$ and 
bijections $T_x:P_U(x)\to\mathcal{T}(x)$ for $x\in X$ 
such that the following conditions hold for each $x\in X$. 
\begin{enumerate}
\item For each $p\in P_U(x)$, 
$B(p,c_2M/(d+1))$ is contained in $T_x(p)$. 
In particular, $\mathcal{T}(x)$ has capacity $c_2M/(d+1)$. 
\item $\mathcal{T}(x)$ is $c_2M/(d+1)$-separated. 
\item For each $p\in P_U(x)$ and $L>0$, 
the $d$-dimensional volume of $B(\partial T_x(p),L)$ is 
less than $10^{d^2}2LV_{d-1}(\sqrt{5}M+L)^{d-1}$. 
\item Each element of $T_x(p)$ meets at most $10^d$ other elements. 
\item Let $0\leq n<d$. 
If $n+1$ distinct elements $t_0,t_1,\dots,t_n$ in $\mathcal{T}(x)$ meet, 
then there exist $t_{n+1},\dots,t_d\in\mathcal{T}(x)$ such that 
$t_0,t_1,\dots,t_d$ are all distinct and 
$t_0\cap t_1\cap\dots\cap t_d$ is not empty. 
\item Let $\tau=\{p_0,p_1,\dots,p_d\}$ be 
$d+1$ distinct points in $P_U(x)$. 
If $T_x(p_0)\cap T_x(p_1)\cap\dots\cap T_x(p_d)$ is not empty, then 
$B(\Conv(\tau),c_2M/(d+1))$ is contained in $\bigcup_iT_x(p_i)$. 
In addition, 
$B(T_x(p_0),c_2M/(d+1))$ does not meet $\Conv(\tau\setminus\{p_0\})$. 
\end{enumerate}
\end{prop}
\begin{proof}
From Proposition \ref{triangle}, 
we get a clopen set $U\subset X$ and 
a family $\{\mathcal{D}(x)\mid x\in X\}$. 
By applying the construction in this section to each $\mathcal{D}(x)$, 
a tessellation $\mathcal{T}(x)=\{T_x(p)\mid p\in P_U(x)\}$ is obtained. 
By the conditions (7) and (8) of Proposition \ref{triangle} 
and the definition of $\mathcal{T}(x)$, 
it is easy to see that 
the family $\mathcal{T}(x)$ for $x\in X$ is $\phi$-regular. 

The conditions (1), (2) and (3) are immediate consequences 
from Lemma \ref{capacity}, \ref{separated} and \ref{boundary}. 
The condition (4) is clear from the proof of Lemma \ref{boundary}. 
The condition (5) is also clear from the proof of Lemma \ref{separated}. 

Let us check the condition (6). 
Let $\tau=\{p_0,p_1,\dots,p_d\}$ be 
$d+1$ distinct points in $P_U(x)$ and 
suppose that 
$T_x(p_0)\cap T_x(p_1)\cap\dots\cap T_x(p_d)$ is not empty. 
From the proof of Lemma \ref{separated}, 
we see that $\tau$ is in $\mathcal{D}(x)$. 
Suppose that $B(\Conv(\tau),c_2M/(d+1))$ meets $T_x(q)$. 
By Lemma \ref{neighbour}, $\Conv(\tau)$ meets $S(q)$. 
It follows that $q$ is in $\tau$. 
Therefore $B(\Conv(\tau),c_2M/(d+1))$ is covered by $T_x(p_i)$'s. 
Since $B(T_x(p_0),c_2M/(d+1))$ is contained 
in the interior of $S(p_0)$ and 
$\Conv(\tau\setminus\{p_0\})$ is in the boundary of $S(p_0)$, 
the last statement also follows easily. 
\end{proof}

\section{Refining tessellations}

In the last section, 
we gave a method of producing $\phi$-regular tessellations. 
The next step is to show how we may produce a sequence 
having larger and larger elements (more and more separated) 
in such a way that 
each element of one is the union of elements from the previous. 
At the same time, we will need several extra technical conditions 
which will be used later in the proof of the main result. 
While we will provide rigorous and fairly complete arguments, 
most of these properties can be seen fairly easily 
by drawing some pictures. 
Moreover, most of our arguments are similar to 
those in \cite[Theorem 5.1]{GMPS1}, 
and so the reader may refer to it. 

Before stating the result, we will need some notation. 
We are considering a tessellation $\mathcal{T}$ of $\R^d$ 
by polyhedral regions with non-overlapping interiors. 
Given a point $p$ in $\R^d$, we would like to say that 
this point belongs to a unique element of $\mathcal{T}$. 
Of course, this is false since the elements overlap on their boundaries. 
To resolve this difficulty in an arbitrary, but consistent way, 
we define, for any $t$ in $\mathcal{T}$, 
\[
t^*=\{p\in\R^d\mid p+(\ep,\ep^2,\dots,\ep^d)\in t
\text{ for all sufficiently small }\ep>0\}, 
\]
where $\ep^i$ is the $i$-th power of $\ep$. 
It is not so hard to see the following. 
\begin{itemize}
\item $t$ contains $t^*$ and $t^*$ contains the interior of $t$. 
\item If $\mathcal{T}$ is a tessellation by polyhedral regions, then 
the collection $\{t^*\mid t\in\mathcal{T}\}$ is a partition of $\R^d$. 
\item $(t+p)^*=t^*+p$ for any $p\in\R^d$. 
\end{itemize}
We remark that any other definition of $t^*$ would work, 
as long as the above properties are satisfied. 

A comment is in order regarding tessellations. 
In the process we are about to undertake, 
we will take unions of polyhedra, which may be disconnected. 
In addition, a vertex in some polyhedron may only belong to 
one other element of the tessellation. 
So we do not use the terms `vertex', `edge' and `face'. 
Instead, we would like to regard these objects in a combinatorial way 
as the (non-empty) intersection of several polyhedra. 

Let $X\subset\Omega$ be a flat Cantor transversal 
of a free minimal action $\phi$ of $\R^d$ on $\Omega$ 
and let $R_\phi\subset X\times X$ be the \'etale equivalence relation 
induced from $(\Omega,\phi)$. 
We suppose that $P_X$ is $M_0$-separated and $M_1$-syndetic. 

\begin{thm}\label{refine}
There exist a sequence of clopen subsets $U_0,U_1,U_2,\dots$ of $X$ 
and a sequence of $\phi$-regular tessellations 
$\mathcal{T}_0,\mathcal{T}_1,\mathcal{T}_2,\dots$ 
satisfying each of the following conditions for all $l\geq0$ and $x$ in $X$. 
\begin{enumerate}
\item $\mathcal{T}_l(x)$ has capacity $(10^d+l+1)M_1$. 
\item $\mathcal{T}_{l+1}(x)$ is $\diam(\mathcal{T}_l)$-separated. 
\item Each element of $\mathcal{T}_l(x)$ meets at most $10^d$ other elements. 
\item If $t_1,t_2,\dots,t_n\in\mathcal{T}_l(x)$ are mutually distinct 
and $t_1\cap t_2\cap\dots\cap t_n$ is non-empty, then 
there exist $t_{n+1},\dots,t_{d+1}$ in $\mathcal{T}_l(x)$ such that 
$t_i\neq t_j$ for $i\neq j$ and 
$t_1\cap t_2\cap\dots\cap t_{d+1}$ is non-empty. 
\item Each element of $\mathcal{T}_l(x)$ is contained in 
an element of $\mathcal{T}_{l+1}(x)$. 
\item For any $s\in\mathcal{T}_{l+1}(x)$, 
$\#\{t\in\mathcal{T}_l(x)\mid t'\subset s
\text{ for any }t'\in\mathcal{T}_l(x)\text{ with }t\cap t'\neq\emptyset\}$ 
is not less than 
$10^{d^2}\#(P_X(x)\cap B(\partial s,\diam(\mathcal{T}_0)))$. 
\item There exists a bijection $\pi_{l,x}:\mathcal{T}_l(x)\to P_{U_l}(x)$ 
such that $\pi_{l,x}(t)$ is in the interior of $t$ 
for each $t\in\mathcal{T}_l(x)$. 
\end{enumerate}
\end{thm}
\begin{proof}
By Proposition \ref{wellsepa}, 
we may choose a clopen subset $U_0$ and 
a $\phi$-regular family of tessellations $\mathcal{T}_0$ 
such that $\diam(\mathcal{T}_0(x))$ has capacity $(10^d+1)M_1$ 
for each $x\in X$. 
The other properties required for $\mathcal{T}_0$ easily follow 
from Proposition \ref{wellsepa}. 

Next, we suppose that we have found a clopen subset $U_l$ and 
a $\phi$-regular family of tessellations $\mathcal{T}_l$ 
satisfying the desired conditions for some $l\geq0$. 
There exists a constant $K>0$ such that $P_{U_l}$ is $K$-syndetic, 
because $U_l$ is clopen. 
Put 
\[
L=\diam(\mathcal{T}_0)+M_0/2+\diam(\mathcal{T}_l). 
\]
Let $\pi_{l,x}:\mathcal{T}_l(x)\to P_{U_l}(x)$ be the bijection 
described in the condition (7). 

We would like to construct $U_{l+1}$ and $\mathcal{T}_{l+1}$. 
We find a constant $M>0$ satisfying each of the following: 
\begin{equation}
c_2(d+1)^{-1}M\geq(10^d+l+2)M_1+\diam(\mathcal{T}_l), 
\label{ineq1}
\end{equation}
\begin{equation}
c_2(d+1)^{-1}M\geq2\diam(\mathcal{T}_l) 
\label{ineq2}
\end{equation}
and 
\begin{equation}
\left(\frac{c_2(d+1)^{-1}M-2\diam(\mathcal{T}_l)-K}{K}\right)^d
\geq10^{d^2}\times
\frac{10^{d^2}2LV_{d-1}(\sqrt{5}M+L)^{d-1}}{V_d(M_0/2)^d}, 
\label{ineq3}
\end{equation}
where $V_k$ denotes the $k$-dimensional volume of 
the unit ball in $\R^k$. 
We note that in the third inequality, 
the left hand side is a polynomial of degree $d$ with variable $M$ and 
the right hand side is a polynomial of degree $d-1$ with variable $M$. 

By applying Proposition \ref{wellsepa} to the constant $M$, 
we may find a clopen subset $U\subset X$, 
a $\phi$-regular family of tessellations $\mathcal{T}$ and 
bijections $\pi_x:\mathcal{T}(x)\to P_U(x)$ for $x$ in $X$ 
satisfying the conditions given there. 
For each $x$ in $X$ and $s$ in $\mathcal{T}(x)$, we define 
\[
\tilde{s}=\bigcup t, 
\]
where the union runs over all the cells $t\in\mathcal{T}_l(x)$ 
such that $\pi_{l,x}(t)$ is in $s^*$. 
We now define a new tessellation $\mathcal{T}_{l+1}(x)$ 
to be the collection of all $\tilde{s}$, where $s$ is in $\mathcal{T}(x)$. 
Since $\mathcal{T}_l$ and $P_{U_l}$ are $\phi$-regular, 
this new family of tessellations $\mathcal{T}_{l+1}$ is 
also $\phi$-regular. 
Clearly, condition (5) is satisfied. 

We first observe that, for any $s$ in $\mathcal{T}(x)$, 
each point in $\tilde{s}$ is in some element $t$ of $\mathcal{T}_l(x)$ 
with $\pi_{l,x}(t)\in s^*$. 
As the diameter of $t$ is at most $\diam(\mathcal{T}_l)$, 
it follows that 
every point of $\tilde{s}$ is within distance $\diam(\mathcal{T}_l)$ of $s$. 

We next verify that $\mathcal{T}_{l+1}(x)$ has capacity $(10^d+l+2)M_1$. 
Each element $s$ in $\mathcal{T}(x)$ has capacity $c_2(d+1)^{-1}M$ 
and hence contains an open ball $B(p,c_2(d+1)^{-1}M)$. 
It follows that 
the open ball $B(p,c_2(d+1)^{-1}M-\diam(\mathcal{T}_l))$ is 
contained in $\tilde{s}$. 
Therefore $\tilde{s}$ has capacity $(10^d+l+2)M_1$ by \eqref{ineq1}. 

We will show that 
the map sending $s$ in $\mathcal{T}(x)$ 
to $\tilde{s}$ in $\mathcal{T}_{l+1}(x)$ 
is a bijection which preserves non-trivial (multiple) intersections. 
The first step in this is to observe that 
if $\tilde{s}_1,\tilde{s}_2,\dots,\tilde{s}_n$ have 
a non-trivial intersection, then so do $s_1,s_2,\dots,s_n$. 
Let $p$ be in the intersection of $\tilde{s}_1,\tilde{s}_2,\dots,\tilde{s}_n$. 
Since every point of $\tilde{s}_i$ is 
within distance $\diam(\mathcal{T}_l)$ of $s_i$, 
we obtain $d(p,s_i)\leq\diam(\mathcal{T}_l)$ for each $i=1,2,\dots,n$. 
By Proposition \ref{wellsepa} (2), 
$\mathcal{T}(x)$ is $c_2(d+1)^{-1}M$-separated. 
It follows from \eqref{ineq2} that 
$s_1,s_2,\dots,s_n$ have a non-trivial intersection. 
In particular, any distinct $d+2$ elements in $\mathcal{T}_{l+1}(x)$ 
do not meet. 

Now, we want to consider the situation that 
$s_1,s_2,\dots,s_n$ are $n$ distinct elements of $\mathcal{T}(x)$ 
with a non-trivial intersection. 
We will show that $\tilde{s}_1,\tilde{s}_2,\dots,\tilde{s}_n$ also 
have a non-trivial intersection. 
From Proposition \ref{wellsepa} (5), 
we may find $s_{n+1},\dots,s_{d+1}$ in $\mathcal{T}(x)$ such that 
$s_1,s_2,\dots,s_{d+1}$ are all distinct and have a non-trivial intersection. 
It suffices to show that 
$\tilde{s}_1,\tilde{s}_2,\dots,\tilde{s}_{d+1}$ have 
a non-trivial intersection. 
The proof is by contradiction. 
We assume that they do not meet. 
For each $i$, we put $p_i=\pi_x(s_i)\in P_U(x)$. 
Let $\tau=\{p_1,p_2,\dots,p_{d+1}\}$. 
From the observation above, for each $i$, 
$\tilde{s_i}$ is contained in $B(s_i,\diam(\mathcal{T}_l))$. 
By Proposition \ref{wellsepa} (6), 
$B(s_i,c_2(d+1)^{-1}M)$ does not meet $\Conv(\tau\setminus\{p_i\})$. 
It follows from \eqref{ineq2} that 
$\tilde{s_i}$ does not meet $\Conv(\tau\setminus\{p_i\})$. 
Besides, by Proposition \ref{wellsepa} (6), 
$B(\Conv(\tau),c_2(d+1)^{-1}M)$ is covered by $s_i$'s. 
If $\Conv(\tau)$ meets $\tilde{s}\in\mathcal{T}_{l+1}(x)$ 
which is distinct from $\tilde{s}_i$'s, then 
$B(\Conv(\tau),c_2(d+1)^{-1}M)$ meets $s$, which is a contradiction. 
Hence $\Conv(\tau)$ is covered by $\tilde{s}_i$'s. 
We define a continuous function $f:\Conv(\tau)\to\R$ by 
\[
f(q)=\sum_{i=1}^{d+1}d(q,\tilde{s}_i). 
\]
Since we have assumed that $\tilde{s}_i$'s have trivial intersections, 
$f$ is strictly positive. 
Define a continuous map $g:\Conv(\tau)\to\Conv(\tau)$ by 
\[
g(q)=f(q)^{-1}\sum_{i=1}^{d+1}d(q,\tilde{s_i})p_i
\]
for $q\in\Conv(\tau)$. 
It is easy to see that $g(q)$ is in $\Conv(\tau)$. 
By the Brouwer fixed point theorem, 
we can find a fixed point $q$ of the map $g$. 
There exists $i$ such that $q$ is in $\tilde{s}_i$, 
because $\Conv(\tau)$ is contained in the union of $\tilde{s}_i$'s. 
By the definition of $g$, one has $g(q)\in\Conv(\tau\setminus\{p_i\})$. 
But, we have observed that 
$\tilde{s_i}$ does not meet $\Conv(\tau\setminus\{p_i\})$. 
This contradiction establishes the desired result. 

Therefore, by conditions (4) and (5) of Proposition \ref{wellsepa}, 
we obtain conditions (3) and (4) of Theorem \ref{refine}. 

We next consider condition (2). 
Let $\tilde{s}_1,\tilde{s}_2,\dots,\tilde{s}_n$ be elements 
in $\mathcal{T}_{l+1}(x)$ and let $p$ be in $\R^d$. 
Suppose that $d(p,\tilde{s}_i)$ is less than $\diam(\mathcal{T}_l)$ 
for each $i=1,2,\dots,n$. 
Since every point of $\tilde{s}_i$ is 
within distance $\diam(\mathcal{T}_l)$ of $s_i$, 
we get $d(p,s_i)<2\diam(\mathcal{T}_l)$. 
From Proposition \ref{wellsepa} (2), 
$\mathcal{T}(x)$ is $c_2(d+1)^{-1}M$-separated. 
It follows from \eqref{ineq2} that 
$s_1,s_2,\dots,s_n$ have a non-trivial intersection. 
Hence $\tilde{s}_1,\tilde{s}_2,\dots,\tilde{s}_n$ also have 
a non-trivial intersection. 

We now consider condition (6). 
Let $\tilde{s}$ be an element in $\mathcal{T}_{l+1}(x)$. 
From Proposition \ref{wellsepa} (1), 
$s$ contains an open ball $B(p,c_2(d+1)^{-1}M)$. 
We first claim 
\begin{align*}
&\{t\in\mathcal{T}_l(x)\mid 
d(\pi_{l,x}(t),p)<c_2(d+1)^{-1}M-2\diam(\mathcal{T}_l)\} \\
&\subset\{t\in\mathcal{T}_l(x)\mid t'\subset\tilde{s}
\text{ for any }t'\in\mathcal{T}_l(x)\text{ with }t\cap t'\neq\emptyset\}. 
\end{align*}
Indeed, if $d(\pi_{l,x}(t),p)$ is less than 
$c_2(d+1)^{-1}M-2\diam(\mathcal{T}_l)$, then 
$d(\pi_{l,x}(t'),p)$ is less than $c_2(d+1)^{-1}M$ 
for any $t'\in\mathcal{T}_l(x)$ such that $t\cap t'\neq\emptyset$. 
It follows that $t'$ is contained in $\tilde{s}$. 
Since $P_{U_l}(x)$ is $K$-syndetic, 
by Lemma \ref{sepasyn} (2), one has 
\begin{align*}
&\#\{t\in\mathcal{T}_l(x)\mid t'\subset\tilde{s}
\text{ for any }t'\in\mathcal{T}_l(x)\text{ with }t\cap t'\neq\emptyset\} \\
&\geq\#(P_{U_l}(x)\cap B(p,c_2(d+1)^{-1}M-2\diam(\mathcal{T}_l))) \\
&>\left(\frac{c_2(d+1)^{-1}M-2\diam(\mathcal{T}_l)-K}{K}\right)^d. 
\end{align*}
We next consider 
\[
E=P_X(x)\cap B(\partial\tilde{s},\diam(\mathcal{T}_0)). 
\]
If $q$ is in $\partial\tilde{s}$, then 
there exists $s_0\in\mathcal{T}(x)$ which is distinct from $s$ 
and $q\in\tilde{s}\cap\tilde{s}_0$. 
Therefore there exist $t,t_0\in\mathcal{T}_l(x)$ such that 
$q\in t\cap t_0$, $\pi_{l,x}(t)\in s^*$ and $\pi_{l,x}(t_0)\in s_0^*$. 
If $q$ is in the interior of $s$, then 
the line segment from $q$ to $\pi_{l,x}(t_0)$ meets the boundary of $s$. 
If $q$ is not in $s$, then 
the line segment from $q$ to $\pi_{l,x}(t)$ meets the boundary of $s$. 
It follows that 
$d(q,\partial s)$ is not greater than $\diam(\mathcal{T}_l)$. 
From this, for any $u\in E$, we get 
\begin{align*}
B(u,M_0/2)&\subset B(\partial\tilde{s},\diam(\mathcal{T}_0)+M_0/2) \\
&\subset B(\partial s,\diam(\mathcal{T}_0)+M_0/2+\diam(\mathcal{T}_l)) \\
&=B(\partial s,L). 
\end{align*}
For any $u\neq u'\in E$, 
$B(u,M_0/2)$ and $B(u',M_0/2)$ are disjoint, 
because $P_X(x)$ is $M_0$-separated. 
This, together with Proposition \ref{wellsepa} (3), implies that 
\[
\#E\times V_d(M_0/2)^d\leq10^{d^2}2LV_{d-1}(\sqrt{5}M+L)^{d-1}. 
\]
It follows from these estimates and \eqref{ineq3} that 
\[
\#\{t\in\mathcal{T}_l(x)\mid t'\subset\tilde{s}
\text{ for any }t'\in\mathcal{T}_l(x)\text{ with }t\cap t'\neq\emptyset\}
\geq10^{d^2}\times\#E, 
\]
which completes the proof of (6). 

As a final point, we put $U_{l+1}=U$ and 
define the map $\pi_{l+1,x}(\tilde{s})=\pi_x(s)$ 
for all $x$ in $X$ and $s\in\mathcal{T}(x)$. 
Let us check (7). 
From Proposition \ref{wellsepa} (1), 
the open ball centred at $\pi_x(s)$ with radius $c_2(d+1)^{-1}M$ 
is contained in $s$. 
Since $c_2(d+1)^{-1}M$ is greater than $\diam(\mathcal{T}_l)$, 
$\pi_x(s)$ is in the interior of $\tilde{s}$ as required. 
\end{proof}

\section{AF equivalence relations and boundaries}

Let $X\subset\Omega$ be a flat Cantor transversal 
of a free minimal action $\phi$ of $\R^d$ on $\Omega$ 
and let $R_\phi\subset X\times X$ be the \'etale equivalence relation 
induced from $(\Omega,\phi)$. 
We suppose that $P_X$ is $M_0$-separated and $M_1$-syndetic. 

In this section, we would like to use our earlier construction of 
a nested sequence of $\phi$-regular tessellations 
to construct the data necessary in the application of the absorption theorem 
to give a proof of the main result. 
This needs, first of all, an AF relation. 
The obvious choice is by using the interiors of the cells in the tessellation. 
These equivalence relations are actually too large. We will refine them 
by considering all $(d+1)$-tuples $t_1,t_2,\dots,t_d,t_{d+1}$ 
which have non-trivial intersection. 
At the same time, 
we also keep track of boundary sets $B^{n+1}$ 
along the $(d-n)$-dimensional faces for each $n=1,2,\dots,d$. 
We will apply the absorption theorem \cite[Theorem 3.2]{M3} $d$ times: 
the $n$-th application of the absorption theorem enlarges 
the equivalence relation along the $(d-n)$-dimensional faces 
for $n=1,2,\dots,d$. 

We need some notation. 
Let $\{\mathcal{T}(x)\mid x\in X\}$ be a $\phi$-regular tessellation. 
For $x$ in $X$ and any $t$ in $\mathcal{T}(x)$, 
we let $N(x,t)$ denote the set of all $t'$ in $\mathcal{T}(x)$, 
including $t$, which intersect $t$. 
For $x$ in $X$ and 
a $(d+1)$-tuple $\xi=(t_1,t_2,\dots,t_{d+1})$ in $\mathcal{T}(x)$, 
the $i$-th coordinate of $\xi$ is denoted by $\xi(i)$. 
For $n=2,3,\dots,d+1$, 
let $S_n$ denote the permutation group on $\{1,2,\dots,n\}$. 
We regard $S_n$ as a subgroup of $S_{d+1}$ in a obvious fashion. 
For $\sigma\in S_{d+1}$ and $\xi=(t_1,t_2,\dots,t_{d+1})$, 
we define $\sigma(\xi)$ by $\sigma(\xi)(i)=\xi(\sigma^{-1}(i))$. 

We begin with our refining sequence of 
$\phi$-regular tessellations $\mathcal{T}_0,\mathcal{T}_1,\mathcal{T}_2,\dots$ 
provided by Theorem \ref{refine} of the last section. 
For any $x$ in $X$, 
we let 
\[
i(x,\cdot):\mathcal{T}_l(x)\to\mathcal{T}_{l+1}(x)
\]
be the unique function such that $i(x,t)\supset t$ 
for every $t\in\mathcal{T}_l(x)$. 
If $k\geq1$, we let $i^k$ denote the composition of 
$k$ functions $i$ mapping $\mathcal{T}_l$ to $\mathcal{T}_{l+k}$ 
for any $l\geq0$. 
By the $\phi$-regularity, 
for any $x\in X$, $l\geq0$, $k\geq1$ and $t\in\mathcal{T}_l(x)$, 
there exists an open neighbourhood $U$ of $x$ such that 
$t\in\mathcal{T}_l(y)$ and $i^k(y,t)=i^k(x,t)$ 
for every $y\in U$. 

For each $x\in X$ and $l\in\N$, 
we let $\mathcal{T}_l^{(1)}(x)$ denote the set of all $(d+1)$-tuples 
$\xi=(t_1,t_2,\dots,t_{d+1})$ in $\mathcal{T}_l(x)$ such that 
$\bigcap_{i=1}^{d+1}t_i$ is non-empty and $t_i\neq t_j$ for $i\neq j$. 
Further, for each $n=2,3,\dots,d+1$, we let 
\[
\mathcal{T}_l^{(n)}(x)=\left\{\xi\in\mathcal{T}_l^{(1)}(x)\mid
\#\{i^k(x,\xi(1)),\dots,i^k(x,\xi(n))\}=n
\text{ for all }k\geq1\right\}. 
\]
In other words, 
$\xi$ is in $\mathcal{T}_l^{(n)}(x)$ if and only if 
$i^k(x,\xi(1)),i^k(x,\xi(2)),\dots,i^k(x,\xi(n))$ are all distinct 
for any $k\geq1$. 
Clearly we have 
\[
\mathcal{T}_l^{(1)}(x)\supset\mathcal{T}_l^{(2)}(x)\supset
\dots\supset\mathcal{T}_l^{(d+1)}(x). 
\]
\bigskip

First, we would like to define a surjective map 
$\theta_x:P_X(x)\to\mathcal{T}_0^{(1)}(x)$ for each $x\in X$ 
in a $\phi$-regular fashion as follows. 
Consider all possible triples 
\[
(t,t^*\cap P_X(x),N(x,t)), 
\]
where $x$ is in $X$ and $t$ is in $\mathcal{T}_0(x)$. 
We consider $(t_1,t_1^*\cap P_X(x_1),N(x_1,t_1))$ and
$(t_2,t_2^*\cap P_X(x_2),N(x_2,t_2))$ to be equivalent 
if they are translates of one another, 
namely that there exists $p\in\R^d$ 
such that 
\[
t_2=t_1+p, \ t_2^*\cap P_X(x_2)=t_1^*\cap P_X(x_1)+p
\text{ and }N(x_2,t_2)=N(x_1,t_1)+p. 
\]
Since $\mathcal{T}_0$ is $\phi$-regular, 
there exists a finite number of equivalence classes. 
We let $\mathcal{P}$ be a finite set 
containing exactly one representative of each equivalence class. 

Let $(t,F,N)$ be in $\mathcal{P}$. 
For $p\in F$, we let $\theta(t,F,N,p)$ be 
a $(d+1)$-tuple $\xi=(t,t_1,t_2,\dots,t_d)$ such that 
$t_i\in N$, $t\cap t_1\cap\dots\cap t_d\neq\emptyset$ and $\#\xi=d+1$. 
We show that we can choose $\theta(t,F,N,p)$ so that 
$\theta(t,F,N,\cdot)$ becomes a surjection 
from $F$ to such $(d+1)$-tuples. 
From Theorem \ref{refine} (3), 
the number of such $(d+1)$-tuples $\xi$ is less than $(10^d)^d$. 
On the other hand, 
$\mathcal{T}_0(x)$ has capacity $(10^d+1)M_1$ by Theorem \ref{refine} (1), 
and so each element of $\mathcal{T}_0(x)$ contains 
a ball of radius $(10^d+1)M_1$. 
It follows from Lemma \ref{sepasyn} (2) that 
the cardinality of $F$ is not less than $(10^d)^d$. 
Hence, we can choose $\theta(t,F,N,p)$ 
so that $\theta(t,F,N,\cdot)$ becomes a surjection 
from $F$ to the $(d+1)$-tuples as above. 

Haven chosen these items for our representative patterns $\mathcal{P}$, 
we extend their definition by translation as follows. 
Let $x$ be in $X$ and $p$ be in $P_X(x)$. 
Take $t\in\mathcal{T}_0(x)$ such that $p\in t^*$. 
We find the unique $q$ in $\R^d$ 
with $(t+q,t^*\cap P_X(x)+q,N(x,t)+q)$ in $\mathcal{P}$ and define 
\[
\theta_x(p)=\theta(t+q,t^*\cap P_X(x)+q,N(x,t)+q,p+q)-q. 
\]
The following lemma follows at once from the definitions 
and we omit the proof. 

\begin{lem}\label{theta}
For any $x\in X$, we have the following. 
\begin{enumerate}
\item $\theta_x:P_X(x)\to\mathcal{T}_0^{(1)}(x)$ is surjective. 
\item For any $p\in P_X(x)$ and $t=\theta_x(p)(1)$, 
one has $p\in t^*$. 
\item For any $p,q\in P_X(x)$, 
$\theta_{\phi^q(x)}(p-q)=\theta_x(p)-q$. 
\item For any $p\in P_X(x)$, 
there exists an open neighbourhood $U\subset X$ of $x$ 
such that $p\in P_X(y)$ and $\theta_y(p)=\theta_x(p)$ 
for all $y\in U$. 
\end{enumerate}
\end{lem}

For each $n=1,2,3,\dots,d+1$, 
we define a subset $B^n$ of $X$ by 
\[
B^n=\{\phi^p(x)\in X\mid
x\in X,p\in P_X(x),\theta_x(p)\in\mathcal{T}_0^{(n)}(x)\}. 
\]
It is easy to see 
\[
X=B^1\supset B^2\supset B^3\supset\dots\supset B^{d+1}. 
\]
The subsets $B^n$ are $d$-dimensional counterparts 
of minimal and maximal paths 
in Bratteli-Vershik models for minimal $\Z$-actions (\cite{HPS}). 
Actually, when $d=1$, one can construct the refining sequence of tessellations 
so that $B^2$ consists only of two points, 
namely the minimal path and the maximal path. 
For general $d$, the `boundary' has a hierarchic structure: 
$B^n$ corresponds to $(d{-}n{+}1)$-dimensional faces of polytopes. 

\begin{lem}\label{B}
For every $n=1,2,3,\dots,d+1$, 
$B^n$ is a non-empty closed subset of $X$. 
\end{lem}
\begin{proof}
Fix $1\leq n\leq d+1$. 
For $x\in X$ and $l\in\N$, we let 
\[
\mathcal{T}_{0,l}^{(n)}(x)=\left\{\xi\in\mathcal{T}_0^{(1)}(x)\mid
\#\{i^l(x,\xi(1)),\dots,i^l(x,\xi(n))\}=n\right\}. 
\]
Then $\mathcal{T}_{0,l}^{(n)}(x)\supset\mathcal{T}_{0,l+1}^{(n)}(x)$ 
and $\mathcal{T}_{0}^{(n)}(x)=\bigcap_l\mathcal{T}_{0,l}^{(n)}(x)$. 
Define 
\[
B_l^n=\{\phi^p(x)\in X\mid
x\in X,p\in P_X(x),\theta_x(p)\in\mathcal{T}_{0,l}^{(n)}(x)\}. 
\]
It is clear that $B_l^n\supset B_{l+1}^n$ and $B^n=\bigcap_lB_l^n$. 
Thus, it suffices to show that $B_l^n$ is non-empty and closed. 

By Theorem \ref{refine} (4), 
there exist $t_1,t_2,\dots,t_{d+1}$ in $\mathcal{T}_l(x)$ such that 
$t_i\neq t_j$ for $i\neq j$ and 
$t_1\cap t_2\cap\dots\cap t_{d+1}\neq\emptyset$. 
Since each element of $\mathcal{T}_l(x)$ is 
a union of elements of $\mathcal{T}_0(x)$, 
we may find $s_1,s_2,\dots,s_{d+1}$ in $\mathcal{T}_0(x)$ such that 
$i^l(x,s_j)=t_j$ and $s_1\cap s_2\cap\dots\cap s_{d+1}\neq\emptyset$. 
It follows that 
$(s_1,s_2,\dots,s_{d+1})$ belongs to $\mathcal{T}_{0,l}^{(d+1)}(x)$. 
Hence $\mathcal{T}_{0,l}^{(n)}(x)$ is non-empty 
for any $n=1,2,\dots,d+1$. 
From this, we see that $B_l^n$ is non-empty. 

Let us show that $B_l^n$ is closed. 
Take $x=\phi^0(x)\in X\setminus B_l^n$ arbitrarily. 
Put $\xi=\theta_x(0)$. 
Since $\xi$ is not in $\mathcal{T}_{0,l}^{(n)}(x)$, 
\[
\#\{i^l(x,\xi(1)),i^l(x,\xi(2)),\dots,i^l(x,\xi(n))\}<n. 
\]
By the $\phi$-regularity, 
there exists an open neighbourhood $U$ of $x$ such that, 
for any $y\in U$, one has 
$\theta_y(0)=\theta_x(0)=\xi$ 
and $i^l(y,\xi(j))=i^l(x,\xi(j))$ for $j=1,2,\dots,n$. 
It follows that 
\[
\#\{i^l(y,\xi(1)),i^l(y,\xi(2)),\dots,i^l(y,\xi(n))\}<n. 
\]
Therefore one has $\theta_y(0)=\xi\notin\mathcal{T}_{0,l}^{(n)}(y)$ 
and $y=\phi^0(y)$ is not in $B_l^n$. 
This then implies that 
$U$ is contained in $X\setminus B_l^n$, and so $B_l^n$ is closed. 
\end{proof}

\bigskip

Next, we would like to introduce equivalence relations $Q^n$ and $R^n$ 
on each $B^n$. 
First, for each $n=2,3,\dots,d+1$, 
we define an equivalence relation $Q^n$ on $B^n$ as follows. 
Let $\phi^p(x)$ and $\phi^q(x)$ be two points in $B^n$, 
where $\theta_x(p)$ and $\theta_x(q)$ are in $\mathcal{T}_0^{(n)}(x)$. 
Put $\xi=\theta_x(p)$ and $\eta=\theta_x(q)$. 
For every $l\in\N$, 
we let the pair $(\phi^p(x),\phi^q(x))$ be in $Q_l^n$, 
if 
\[
\{i^l(x,\xi(k))\mid k=1,2,\dots,n\}
=\{i^l(x,\eta(k))\mid k=1,2,\dots,n\}. 
\]
It is easy to see that 
$Q_l^n$ is an equivalence relation on $B^n$ 
and $Q_l^n\subset Q_{l+1}^n$. 
We define 
\[
Q^n=\bigcup_{l\in\N}Q_l^n. 
\]

\begin{lem}\label{Q=AF}
Equipped with the relative topology from $R_\phi$, 
the equivalence relation $Q^n$ on $B^n$ is an AF relation 
for each $n=2,3,\dots,d+1$. 
\end{lem}
\begin{proof}
It suffices to show that $Q_l^n$ is 
a compact \'etale relation on $B^n$ for each $l\in\N$. 
First, let us verify the \'etaleness of $Q_l^n$. 
Take $(\phi^p(x),\phi^q(x))\in Q_l^n$, 
where $\xi=\theta_x(p)$ and $\eta=\theta_x(q)$ 
are in $\mathcal{T}_0^{(n)}(x)$ and 
\[
\{i^l(x,\xi(k))\mid k=1,2,\dots,n\}
=\{i^l(x,\eta(k))\mid k=1,2,\dots,n\}. 
\]
There exists a clopen neighbourhood $U$ of $x$ in $X$ 
such that $p,q\in P_X(y)$ for all $y\in U$. 
Put 
\[
O=\{(\phi^p(y),\phi^q(y))\in R_\phi\mid y\in U\}. 
\]
Then $O$ is a clopen neighbourhood of $(\phi^p(x),\phi^q(x))$ 
in $R_\phi$ and the maps 
\[
r:(\phi^p(y),\phi^q(y))\mapsto\phi^p(y)
\quad\text{and}\quad
s:(\phi^p(y),\phi^q(y))\mapsto\phi^q(y)
\]
are local homeomorphisms 
from $O$ to $\phi^p(U)$ and $\phi^q(U)$, respectively. 
From Lemma \ref{theta} (4) and the definition of $i$, 
we may also assume that the neighbourhood $U$ of $x$ in $X$ 
is chosen sufficiently small so that, for any $y\in U$,  
one has $\xi=\theta_y(p)$, $\eta=\theta_y(q)$ and 
$i^l(y,\xi(k))=i^l(x,\xi(k))$,  $i^l(y,\eta(k))=i^l(x,\eta(k))$ 
for each $k=1,2,\dots,n$. 

It is clear that 
$O\cap Q_l^n$ is a clopen neighbourhood of $(\phi^p(x),\phi^q(x))$ 
in $Q_l^n$. 
In order to prove the \'etaleness of $Q_l^n$, 
we would like to show that 
the restriction of the map $r$ to $O\cap Q_l^n$ is 
a homeomorphism to $\phi^p(U)\cap B^n$, and 
it suffices to show that this map is a surjection. 
Choose $y\in U$ and suppose that $\phi^p(y)$ is in $B^n$. 
It follows that $\theta_y(p)=\xi$ is in $\mathcal{T}_0^{(n)}(y)$, 
and so we have 
\[
\#\{i^m(y,\xi(1)),i^m(y,\xi(2)),\dots,i^m(y,\xi(n))\}=n
\]
for every $m\geq0$. 
Since 
\begin{align*}
& \{i^l(y,\xi(1)),i^l(y,\xi(2)),\dots,i^l(y,\xi(n))\} \\
&=\{i^l(x,\xi(1)),i^l(x,\xi(2)),\dots,i^l(x,\xi(n))\} \\
&=\{i^l(x,\eta(1)),i^l(x,\eta(2)),\dots,i^l(x,\eta(n))\} \\
&=\{i^l(y,\eta(1)),i^l(y,\eta(2)),\dots,i^l(y,\eta(n))\}, 
\end{align*}
we can see that 
\[
\#\{i^m(y,\eta(1)),i^m(y,\eta(2)),\dots,i^m(y,\eta(n))\}=n
\]
for every $m\geq0$, that is, $\phi^q(y)$ is in $B^n$. 
In addition, we have $(\phi^p(y),\phi^q(y))$ is in $Q_l^n$. 
Thus, $(\phi^p(y),\phi^q(y))$ is in $O\cap Q_l^n$, 
and so the map $r$ is a local homeomorphism and 
$Q_l^n$ is \'etale. 

Similarly, it can be easily shown that 
$Q_l^n$ is a closed subset of $R_\phi$. 

We next verify that $Q_l^n$ is compact. 
Take $(\phi^p(x),\phi^q(x))\in Q_l^n$, 
where $\xi=\theta_x(p)$ and $\eta=\theta_x(q)$ 
are in $\mathcal{T}_0^{(n)}(x)$ and 
\[
\{i^l(x,\xi(k))\mid k=1,2,\dots,n\}
=\{i^l(x,\eta(k))\mid k=1,2,\dots,n\}. 
\]
From Lemma \ref{theta} (2), 
we have $p\in\xi(1)^*$ and $q\in\eta(1)^*$. 
By definition, one has 
\[
\xi(1)^*\subset\xi(1)\subset i^l(x,\xi(1))
\quad\text{and}\quad
\eta(1)^*\subset\eta(1)\subset i^l(x,\eta(1)). 
\]
Since $i^l(x,\xi(1))$ meets $i^l(x,\eta(1))$, 
we can conclude that $d(p,q)\leq2\diam(\mathcal{T}_l(x))$. 
It follows from the $\phi$-regularity of $\mathcal{T}_l$ that 
$\diam(\mathcal{T}_l(x))$ is bounded uniformly over all $x\in X$. 
Therefore $Q_l^n$ is compact. 
\end{proof}

We need to introduce another equivalence relation $R^n$ 
on each $B^n$ for $n=1,2,\dots,d+1$. 
The relation $R^n$ will be defined by using functions 
\[
j_n(x,\cdot):\mathcal{T}_l^{(n)}(x)\to\mathcal{T}_{l+1}^{(n)}(x)
\]
for every $x\in X$, $l\geq0$ and $n=1,2,\dots,d+1$. 
To define the map $j_n$, we need some notation. 

For each $l\geq0$, 
let $U_l$ be the clopen set corresponding to 
the $\phi$-regular tessellation $\mathcal{T}_l$ 
as described in Theorem \ref{refine}. 
There exists a bijection between $P_{U_l}(x)$ and $\mathcal{T}_l(x)$ 
for each $x$ in $X$. 
As mentioned in Section 4, 
we equip $\R^d$ with the lexicographic ordering. 
By transferring the lexicographic ordering on $P_{U_l}(x)\subset\R^d$, 
we can equip $\mathcal{T}_l(x)$ with a linear ordering. 
By the translation invariance of the lexicographic ordering and 
the $\phi$-regularity of $P_{U_l}$ and $\mathcal{T}_l$, 
we obtain the $\phi$-regularity of the linear ordering 
on $\mathcal{T}_l(x)$ in the following sense. 
\begin{itemize}
\item For any $x\in X$ and $p\in P_X(x)$, 
if $t,t'\in\mathcal{T}_l(x)$ satisfy $t<t'$, 
then $t-p,t'-p\in\mathcal{T}_l(\phi^p(x))$ satisfy $t-p<t'-p$. 
\item For any $x\in X$ and $p\in P_X(x)$, 
there exists an open neighbourhood $U$ of $x$ in $X$ such that, 
if $t,t'\in\mathcal{T}_l(x)$ satisfy $t<t'$, 
then $t,t'\in\mathcal{T}_l(y)$ satisfy $t<t'$ in $\mathcal{T}_l(y)$ 
for any $y\in U$. 
\end{itemize}
By using this linear ordering on $\mathcal{T}_l(x)$, 
we can equip $\mathcal{T}_l^{(1)}(x)$ with the lexicographic ordering. 
It is clear that 
the ordering on $\mathcal{T}_l^{(1)}(x)$ is again $\phi$-regular 
in an obvious sense. 

We would like to define a map $a_l(x,\cdot)$ for each $l\geq0$. 
Fix $l\geq0$. 
For each $x\in X$, 
let $I_l(x)\subset\mathcal{T}_l(x)$ be 
the subset consisting of all $t\in\mathcal{T}_l(x)$ 
such that $i(x,N(x,t))=\{i(x,t)\}$. 
Let $a_l(x,\cdot):I_l(x)\to\mathcal{T}_{l+1}^{(1)}(x)$ 
be a map which satisfies the following. 
\begin{itemize}
\item For any $t\in I_l(x)$, we have $a_l(x,t)(1)=i(x,t)$. 
\item For any $t\in I_l(x)$ and $p\in P_X(x)$, 
we have $a_l(\phi^p(x),t-p)=a_l(x,t)-p$. 
\item For any $t\in I_l(x)$, 
there exists an open neighbourhood $U\subset X$ of $x$ 
such that $t\in I_l(y)$ and $a_l(y,t)=a_l(x,t)$ for every $y\in U$. 
\item For any $\zeta\in\mathcal{T}_{l+1}^{(1)}(x)$ and $s=\zeta(1)$, 
we have 
\[
\#\{t\in I_l(x)\mid a_l(x,t)=\zeta\}\geq
\#(P_X(x)\cap B(\partial s,\diam(\mathcal{T}_0(x)))). 
\]
\end{itemize}
The second and third conditions follow from 
the $\phi$-regularity of $\mathcal{T}_l$ and $\mathcal{T}_{l+1}$. 
For $s\in\mathcal{T}_{l+1}(x)$, from Theorem \ref{refine} (3), 
the number of $(d+1)$-tuples $\zeta\in\mathcal{T}_{l+1}^{(1)}(x)$ 
such that $\zeta(1)=s$ is less than $10^{d^2}$. 
This, together with Theorem \ref{refine} (6), implies 
the fourth condition. 

We now turn to the definition of the map 
$j_n(x,\cdot):\mathcal{T}_l^{(n)}(x)\to\mathcal{T}_{l+1}^{(n)}(x)$. 
Fix $1\leq n\leq d+1$, $l\geq0$ and $x$ in $X$. 
Take $\xi\in\mathcal{T}_l^{(n)}(x)$. 
Put 
\[
M=\{i(x,\xi(1)),i(x,\xi(2)),\dots,i(x,\xi(d+1))\}
\]
and 
\[
N=\{i(x,t)\mid t\in\mathcal{T}_l(x),
t\cap\xi(1)\cap\xi(2)\cap\dots\cap\xi(n)\neq\emptyset\}. 
\]
Clearly $M\subset N$. 
By definition of $\mathcal{T}_l^{(n)}(x)$, one has $n\leq\#M$. 
It follows from Theorem \ref{refine} (2) that 
the elements in $N$ have non-trivial intersection 
and the cardinality of $N$ is not greater than $d+1$. 

First, if $\#N=1$ (this automatically implies $n=\#M=1$), then 
we define $j_1(x,\xi)=a_l(x,\xi(1))$. 

Let us consider the case of $\#N\geq2$. 
For $m=1,2,\dots,\#M$, we let 
\[
k_m=\min\{k\mid\#\{i(x,\xi(1)),\dots,i(x,\xi(k))\}=m\}. 
\]
It is easily verified that $k_1=1,k_2=2,\dots,k_n=n$ and 
$k_n<k_{n+1}<\dots<k_{\#M}\leq d+1$. 
Consider all the $(d+1)$-tuples $\zeta\in\mathcal{T}_{l+1}^{(1)}(x)$ 
such that 
\[
\zeta(m)=i(x,\xi(k_m))\text{ for all }m=1,2,\dots,\#M
\]
and 
\[
\{\zeta(1),\zeta(2),\dots,\zeta(\#N)\}=N. 
\]
Define $j_n(x,\xi)$ to be the minimum element 
in the set of such $(d+1)$-tuples $\zeta$ 
with respect to the linear ordering on $\mathcal{T}_{l+1}^{(1)}(x)$. 
From Theorem \ref{refine} (4), such a $(d+1)$-tuple exists. 
Moreover, from Theorem \ref{refine} (3), 
there exist only finitely many such $(d+1)$-tuples. 
Hence $j_n(x,\xi)\in\mathcal{T}_{l+1}^{(1)}(x)$ is well-defined. 
Note that for any $k>\#N$, $j_n(x,\xi)(k)$ is the minimum element 
in the set of all tiles in $\mathcal{T}_{l+1}(x)$ 
which meet $j_n(x,\xi)(1)\cap j_n(x,\xi)(2)\cap\dots\cap j_n(x,\xi)(k-1)$. 

From the definition above, one has $j_n(x,\xi)(k)=i(x,\xi(k))$ 
for any $\xi\in\mathcal{T}_l^{(n)}(x)$ and $k=1,2,\dots,n$. 
Therefore $j_n(x,\xi)$ belongs to $\mathcal{T}_{l+1}^{(n)}(x)$. 
In addition, it is also easy to see that 
$j_n(x,\cdot)$ is $\phi$-regular in an appropriate sense. 
We remark that $j_n(x,\xi)(k)$ for $k>\#N$ depends only on $N$. 

The following lemma is an easy consequence of the definitions 
and we omit the proof. 

\begin{lem}\label{j}
For $1\leq n\leq d+1$, $l\geq0$ and $x\in X$, 
we have the following. 
\begin{enumerate}
\item For any  $\xi\in\mathcal{T}_l^{(m)}(x)$ and $n\leq m$, 
one has $j_n(x,\xi)(k)=i(x,\xi(k))$ for each $k=1,2,\dots,m$. 
\item For any $\xi\in\mathcal{T}_l^{(n)}(x)$ and $m=1,2,\dots,d+1$, 
one has 
\[
\{j_n(x,\xi)(k)\mid k=1,2,\dots,m'\}
=\{i(x,\xi(k))\mid k=1,2,\dots,m\}, 
\]
where $m'=\#\{i(x,\xi(1)),\dots,i(x,\xi(m))\}$. 
\item If $\xi,\eta\in\mathcal{T}_l^{(n)}(x)$ satisfy 
\[
\{\xi(k)\mid k=1,2,\dots,n\}=\{\eta(k)\mid k=1,2,\dots,n\}, 
\]
then there exists $\sigma\in S_{d+1}$ such that 
\[
\sigma(j_n(x,\xi))=j_n(x,\eta)\text{ and }
\sigma(\{1,2,\dots,n\})=\{1,2,\dots,n\}. 
\]
\item For $\xi\in\mathcal{T}_l^{(n)}(x)$ and $\sigma\in S_m$ 
with $\sigma(\{1,2,\dots,n\})=\{1,2,\dots,n\}$, 
there exists $\sigma'\in S_{m'}$ such that 
\[
\sigma'(j_n(x,\xi))=j_n(x,\sigma(\xi))
\]
and 
\[
\sigma'(k)=\sigma(k)\text{ for all }k=1,2,\dots,n, 
\]
where $m'=\#\{i(x,\xi(1)),\dots,i(x,\xi(m))\}$. 
\end{enumerate}
\end{lem}

For $k\geq1$, 
we let $j_n^k(x,\cdot)$ denote the composition of 
$k$ functions $j_n(x,\cdot)$ 
mapping $\mathcal{T}_l^{(n)}(x)$ to $\mathcal{T}_{l+k}^{(n)}(x)$ 
for any $l\geq0$. 

For each $n=1,2,\dots,d+1$, 
we define an equivalence relation $R^n$ on $B^n$ as follows. 
Let $\phi^p(x)$ and $\phi^q(x)$ be two points in $B^n$, 
where $\theta_x(p)$ and $\theta_x(q)$ are in $\mathcal{T}_0^{(n)}(x)$. 
Put $\xi=\theta_x(p)$ and $\eta=\theta_x(q)$. 
For every $l\in\N$, 
we let the pair $(\phi^p(x),\phi^q(x))$ be in $R_l^n$, 
if there exists $\sigma\in S_n$ such that 
\[
\sigma(j_n^l(x,\xi))=j_n^l(x,\eta). 
\]
It is easy to see that 
$R_l^n$ is an equivalence relation on $B^n$. 
If $\sigma(j_n^l(x,\xi))=j_n^l(x,\eta)$, then 
by applying Lemma \ref{j} (4) to the case of $m{=}n$, 
we get 
\begin{align*}
\sigma(j_n^{l+1}(x,\xi))
&=\sigma(j_n(x,j_n^l(x,\xi)))=j_n(x,\sigma(j_n^l(x,\xi))) \\
&=j_n(x,j_n^l(x,\eta))=j_n^{l+1}(x,\eta). 
\end{align*}
Hence $R_l^n$ is contained in $R_{l+1}^n$. 
Define an equivalence relation $R^n$ on $B^n$ by 
\[
R^n=\bigcup_{l\in\N}R_l^n. 
\]
The following lemma can be shown 
in a similar way to Lemma \ref{Q=AF}. 
We omit the proof. 

\begin{lem}\label{R=AF}
Equipped with the relative topology from $R_\phi$, 
the equivalence relation $R^n$ on $B^n$ is an AF relation 
for each $n=1,2,\dots,d+1$. 
\end{lem}

\bigskip

We wish to establish several facts about $R^n$ and $Q^n$. 
Let us collect notation and terminology 
about equivalence relations. 
Let $R$ be an equivalence relation on $X$. 
For a subset $A\subset X$, 
we set 
\[
R[A]=\{x\in X\mid
\text{there exists }y\in A\text{ such that }(x,y)\in R\}. 
\]
For $x\in X$, we denote $R[\{x\}]$ by $R[x]$ and 
call it the $R$-orbit of $x$. 
For a subset $Y\subset X$, 
we let $R|Y$ denote $R\cap(Y\times Y)$ 
and call it the restriction of $R$ to $Y$. 
Suppose that $R$ is equipped with a topology 
in which $R$ is \'etale (\cite[Definition 2.1]{GPS2}). 
A closed subset $Y\subset X$ is called 
$R$-\'etale or \'etale with respect to $R$, 
if the restriction $R|Y$ with the relative topology from $R$ is \'etale. 
A subset $Y\subset X$ is called $R$-thin, 
if $\mu(Y)$ is zero 
for any $R$-invariant probability measure $\mu$ on $X$. 

We begin with the following. 

\begin{lem}\label{R1minimal}
The AF relation $R^1$ on $B^1=X$ is minimal. 
\end{lem}
\begin{proof}
Take $x\in X$ and a non-empty clopen subset $U\subset X$ 
arbitrarily. 
It suffices to show that $R^1[x]$ meets $U$. 
There exists $M>0$ such that $P_U(x)$ is $M$-syndetic. 

Let $l$ be a natural number such that $(10^d+l+1)M_1\geq M$. 
Put $\zeta=j_1^{l+1}(x,\theta_x(0))$. 
By the definition of $a_l(x,\cdot)$, 
there exists $t\in I_l(x)$ such that $a_l(x,t)=\zeta$. 
It follows from Theorem \ref{refine} (1) that 
the tile $t$ contains an open ball of radius $(10^d+l+1)M_1$. 
Hence $P_U(x)$ meets the interior of $t$. 
Thus, there exists $p\in P_U(x)$ such that 
$i^l(x,\theta_x(p)(1))=t$. 
From this, we get $j_1^l(x,\theta_x(p))(1)=t$, 
which implies $j_1^{l+1}(x,\theta_x(p))=a_l(x,t)=\zeta$. 
Therefore $(x,\phi^p(x))$ is in $R_{l+1}^1$. 
Since $\phi^p(x)$ is in $U$, $R^1[x]$ meets $U$ as required. 
\end{proof}

\begin{lem}\label{B2thin}
The closed subset $B^2\subset X$ is $R^1$-thin. 
\end{lem}
\begin{proof}
For each $l\in\N$, let $B_l^2$ be as in Lemma \ref{B}. 
Since $\bigcap_lB_l^2=B^2$, 
it suffices to show $\mu(B_l^2)\to0$ as $l\to\infty$ 
for any $R^1$-invariant measure $\mu$. 
Fix $l\in\N$ and $x\in X$. 
We will compare $\#R_{l+1}^1[x]$ and $\#(R_{l+1}^1[x]\cap B_{l+1}^2)$. 
Let $\zeta=j_1^{l+1}(x,\theta_x(0))$ and $s=\zeta(1)$. 

Suppose that $t\in\mathcal{T}_l(x)$ is in $I_l(x)$ 
and $a_l(x,t)=\zeta$. 
From the definition of $a_l(x,\cdot)$, 
we notice that 
the number of such tiles $t$ is not less than 
\[
\#(P_X(x)\cap B(\partial s,\diam(\mathcal{T}_0(x)))). 
\]
From the definition of $j_1(x,\cdot)$, 
for any $\xi\in\mathcal{T}_l^{(1)}(x)$ such that $\xi(1)=t$, 
one has $j_1(x,\xi)=\zeta$. 
If $p\in P_X(x)$ is contained in the interior of $t$, 
then $i^l(x,\theta_x(p)(1))$ is equal to $t$ 
by Lemma \ref{theta} (2) and the definition of $i(x,\cdot)$. 
By the repeated use of Lemma \ref{j} (1), 
we get $j_1^l(x,\theta_x(p))(1)=t$. 
Therefore $j_1^{l+1}(x,\theta_x(p))=\zeta$. 
Thus, if $p\in P_X(x)$ is contained in the interior of $t$, 
then $(x,\phi^p(x))$ is in $R_{l+1}^1$. 
By Theorem \ref{refine} (1), 
$t$ contains an open ball of radius $(10^d+l+1)M_1$. 
It follows from Lemma \ref{sepasyn} (2) that 
the number of points in the intersection of 
$P_X(x)$ and the interior of $t$ is not less than $(10^d+l)^d$. 
Hence we have 
\begin{align*}
\#R_{l+1}^1[x]&=\#\{p\in P_X(x)\mid(x,\phi^p(x))\in R_{l+1}^1\} \\
&\geq(10^d+l)^d
\times\#(P_X(x)\cap B(\partial s,\diam(\mathcal{T}_0(x)))). 
\end{align*}

Next, suppose that 
$\phi^p(x)$ is in $R_{l+1}^1[x]\cap B_{l+1}^2$. 
Put $\eta=\theta_x(p)$. 
By Lemma \ref{theta} (2), $p$ belongs to $\eta(1)^*$. 
From $(x,\phi^p(x))\in R_{l+1}^1$, 
one has $j_1^{l+1}(x,\eta)=\zeta$, and so 
\[
i^{l+1}(x,\eta(1))=j_1^{l+1}(x,\eta)(1)=\zeta(1)=s. 
\]
Since $\phi^p(x)$ is in $B_{l+1}^2$, 
we have $i^{l+1}(x,\eta(1))\neq i^{l+1}(x,\eta(2))$. 
By the definition of the map $i$, 
$\eta(1)\subset i^{l+1}(x,\eta(1))$ and 
$\eta(2)\subset i^{l+1}(x,\eta(2))$. 
Therefore 
\[
d(p,\partial s)\leq d(p,s\cap i^{l+1}(x,\eta(2)))
\leq d(p,\eta(1)\cap\eta(2))<\diam(\mathcal{T}_0(x)). 
\]
It follows that 
\[
\#\{p\in P_X(x)\mid\phi^p(x)\in R_{l+1}^1[x]\cap B_{l+1}^2\}
\leq\#(P_X(x)\cap B(\partial s,\diam(\mathcal{T}_0(x)))). 
\]

Consequently, 
\[
\mu(B_{l+1}^2)\leq(10^d+l)^{-d}
\]
for any $R^1$-invariant probability measure $\mu$ on $X$, 
which completes the proof. 
\end{proof}

\begin{lem}\label{Rn<Qn}
For each $n=2,3,\dots,d+1$, 
$R^n$ is an open subrelation of $Q^n$. 
\end{lem}
\begin{proof}
To show $R^n\subset Q^n$, 
let $(\phi^p(x),\phi^q(x))$ be a pair in $R^n$. 
Let $\xi=\theta_x(p)$ and $\eta=\theta_x(q)$. 
There exists a natural number $l$ and $\sigma\in S_n$ such that 
$\sigma(j_n^l(x,\xi))=j_n^l(x,\eta)$. 
By the repeated use of Lemma \ref{j} (1), we get 
\[
j_n^l(x,\xi)(k)=i^l(x,\xi(k))
\quad\text{and}\quad j_n^l(x,\eta)(k)=i^l(x,\eta(k))
\]
for any $k=1,2,\dots,n$. 
Hence 
\[
\{i^l(x,\xi(k))\mid k=1,2,\dots,n\}
=\{i^l(x,\eta(k))\mid k=1,2,\dots,n\}, 
\]
which implies that $(\phi^p(x),\phi^q(x))$ is in $Q^n$. 

Since both $R^n$ and $Q^n$ are \'etale 
with the induced topology from $R_\phi$, 
we see that $R^n$ is open in $Q^n$. 
\end{proof}

\begin{lem}\label{Bn+1etale}
For each $n=1,2,\dots,d$, 
the closed subset $B^{n+1}$ is $R^n$-\'etale. 
\end{lem}
\begin{proof}
Let $(\phi^p(x),\phi^q(x))$ be a pair in $R^n|B^{n+1}$. 
It suffices to show that 
there exists a clopen neighbourhood $U\subset X$ of $x$ such that 
for any $y\in U$, 
if $(\phi^p(y),\phi^q(y))$ is in $R^n$ and $\phi^p(y)$ is in $B^{n+1}$, 
then $\phi^q(y)$ is also in $B^{n+1}$. 

Put $\xi=\theta_x(p)$ and $\eta=\theta_x(q)$. 
Since $(\phi^p(x),\phi^q(x))$ is in $R^n|B^{n+1}$, 
$\xi$ and $\eta$ are in $\mathcal{T}_0^{(n+1)}(x)$ and 
there exists $l\in\N$ and $\sigma\in S_n$ such that 
$\sigma(j_n^l(x,\xi))=j_n^l(x,\eta)$. 
By the repeated use of Lemma \ref{j} (1), we get 
\[
j_n^l(x,\xi)(k)=i^l(x,\xi(k))
\quad\text{and}\quad j_n^l(x,\eta)(k)=i^l(x,\eta(k))
\]
for any $k=1,2,\dots,n+1$. 
It follows that 
\[
\{i^l(x,\xi(k))\mid k=1,2,\dots,n+1\}
=\{i^l(x,\eta(k))\mid k=1,2,\dots,n+1\}. 
\]
From Lemma \ref{theta} (4), 
there exists a clopen neighbourhood $U$ of $x$ such that 
for any $y\in U$, $\xi=\theta_y(p)$ and $\eta=\theta_y(q)$. 
We may also assume that $U$ is chosen sufficiently small so that, 
for any $y\in U$ and $k=1,2,\dots,n+1$, 
\[
i^l(y,\xi(k))=i^l(x,\xi(k))\quad\text{and}\quad
i^l(y,\eta(k))=i^l(x,\eta(k)). 
\]

Suppose that $y$ is in $U$, $(\phi^p(y),\phi^q(y))$ is in $R^n$ 
and $\phi^p(y)$ is in $B^{n+1}$. 
We get 
\begin{align*}
&\{i^l(y,\eta(k))\mid k=1,2,\dots,n+1\} \\
&=\{i^l(x,\eta(k))\mid k=1,2,\dots,n+1\} \\
&=\{i^l(x,\xi(k))\mid k=1,2,\dots,n+1\} \\
&=\{i^l(y,\xi(k))\mid k=1,2,\dots,n+1\}. 
\end{align*}
Combining this with $\phi^p(y)\in B^{n+1}$, 
we can conclude that $\phi^q(y)$ is in $B^{n+1}$, 
which completes the proof. 
\end{proof}

\begin{lem}\label{Rn<Qn+1}
For each $n=1,2,\dots,d$, 
$R^n|B^{n+1}$ is an open subrelation of $Q^{n+1}$. 
\end{lem}
\begin{proof}
From the proof of the lemma above, 
it is clear that $R^n|B^{n+1}$ is contained in $Q^{n+1}$. 
Since both relations are \'etale 
with the induced topology from $R_\phi$, 
we see that $R^n|B^{n+1}$ is open in $Q^{n+1}$. 
\end{proof}

\begin{lem}\label{RnBn+1=Rphi}
For each $n=1,2,\dots,d$, 
we have $R^n[B^{n+1}]=B^n\cap R_\phi[B^{n+1}]$. 
\end{lem}
\begin{proof}
Clearly $R^n[B^{n+1}]$ is a subset of $B^n\cap R_\phi[B^{n+1}]$, 
and so it suffices to show the other inclusion. 

Let $x$ be in $B^n\cap R_\phi[B^{n+1}]$. 
Put $\xi=\theta_x(0)$. 
For any $l\in\N$, we have 
\[
\#\{i^l(x,\xi(k))\mid k=1,2,\dots,n\}=n. 
\]
From Lemma \ref{j} (1), 
we also get $j_n^l(x,\xi)(k)=i^l(x,\xi(k))$ 
for any $l\in\N$ and $k=1,2,\dots,n$. 
There exists $p\in P_X(x)$ such that $\phi^p(x)\in B^{n+1}$. 
Put $\eta=\theta_x(p)$. 
Similarly, we have 
\[
\#\{i^l(x,\eta(k))\mid k=1,2,\dots,n+1\}=n+1
\]
for any $l\in\N$. 
Choose $L\in\N$ sufficiently large so that 
\[
\diam(\mathcal{T}_{L-1})>
\max\{d(0,\xi(k))\mid k=1,2,\dots,n\}
\]
and 
\[
\diam(\mathcal{T}_{L-1})>
\max\{d(0,\eta(k))\mid k=1,2,\dots,n+1\}. 
\]
From Theorem \ref{refine} (2), 
$\mathcal{T}_L(x)$ is $\diam(\mathcal{T}_{L-1})$-separated. 
Hence we see that, for any $k=1,2,\dots,n+1$, 
$i^L(x,\eta(k))$ meets 
\[
i^L(x,\xi(1))\cap i^L(x,\xi(2))\cap\dots\cap i^L(x,\xi(n)). 
\]
Define 
\[
M=\{j_n^{L+1}(x,\xi)(k)\mid k=1,2,\dots,d+1\}. 
\]
It follows from the definition of the map $j_n$ that 
$i^{L+1}(x,\eta(k))$ is in $M$ for each $k=1,2,\dots,n+1$. 
For each $l>L$, we put 
\[
m_l=\#\{i^{l-L}(x,t)\mid t\in M\}. 
\]
It is easily seen that 
$(m_l)_{l>L}$ is a decreasing sequence of positive integers. 
Therefore, the limit exists. 
In addition, from the argument above, 
$m_l$ is not less than $n+1$, and so 
the limit is not less than $n+1$. 
By the repeated use of Lemma \ref{j} (2), 
for all $l>L$, we have 
\[
\{j_n^l(x,\xi)(k)\mid k=1,2,\dots,m_l\}
=\{i^{l-L}(x,t)\mid t\in M\}. 
\]
Hence, there exists $L'>L$ such that 
$j_n^{L'}(x,\xi)$ is in $\mathcal{T}_{L'}^{(n+1)}(x)$. 
Since we can write each element of $\mathcal{T}_{L'}(x)$ 
as a union of elements of $\mathcal{T}_0(x)$, 
there exists $\zeta\in\mathcal{T}_0^{(n+1)}(x)$ such that 
\[
i^{L'}(x,\zeta(k))=j_n^{L'}(x,\xi)(k)
\]
for all $k=1,2,\dots,d+1$. 
From this we have $j_n^{L'}(x,\zeta)=j_n^{L'}(x,\xi)$. 
By Lemma \ref{theta} (1), 
we can find $q\in P_X(x)$ such that $\theta_x(q)=\zeta$. 
It is easy to see that 
$\phi^q(x)$ is in $B^{n+1}$ and $(x,\phi^q(x))$ is in $R^n$. 
\end{proof}

For each $n=1,2,\dots,d$, 
we define a subset $C^n$ of $X$ by 
\[
C^n=R_\phi[B^n]\setminus R_\phi[B^{n+1}]. 
\]
We let $C^{d+1}=R_\phi[B^{d+1}]$. 
Clearly $C^1,C^2,\dots,C^{d+1}$ are 
mutually disjoint and $R_\phi$-invariant. 
Besides, their union is equal to $X$. 

\begin{lem}\label{RnCn=Rphi}
For each $n=1,2,\dots,d+1$, 
we have $R^n|(B^n\cap C^n)=R_\phi|(B^n\cap C^n)$. 
\end{lem}
\begin{proof}
Clearly $R^n|(B^n\cap C^n)$ is contained in $R_\phi|(B^n\cap C^n)$, 
and so it suffices to show the other inclusion. 

Let $(\phi^p(x),\phi^q(x))$ be a pair in $R_\phi|(B^n\cap C^n)$. 
Put $\xi=\theta_x(p)$ and $\eta=\theta_x(q)$. 
From $\phi^p(x),\phi^q(x)\in B^n$, for any $l\in\N$, one has 
\[
\#\{i^l(x,\xi(k))\mid k=1,2,\dots,n\}=n
\]
and 
\[
\#\{i^l(x,\eta(k))\mid k=1,2,\dots,n\}=n. 
\]

We first claim that 
for any finite subset $M\subset\mathcal{T}_0(x)$, 
there exists $L>0$ such that 
\[
\#\{i^L(x,t)\mid t\in M\}\leq n. 
\]
The proof is by contradiction. 
Suppose that 
$M\subset\mathcal{T}_0(x)$ is a finite subset satisfying 
\[
\#\{i^l(x,t)\mid t\in M\}>n
\]
for any $l>0$. 
By taking a subset of $M$ if necessary, 
we may assume that $M=\{t_1,t_2,\dots,t_{n+1}\}$ and 
\[
\#\{i^l(x,t_k)\mid k=1,2,\dots,n+1\}=n+1
\]
for any $l>0$. 
Choose $L\in\N$ sufficiently large so that 
\[
\diam(\mathcal{T}_{L-1})
>\max\{d(0,t_k)\mid k=1,2,\dots,n+1\}. 
\]
From Theorem \ref{refine} (2), 
$\mathcal{T}_L(x)$ is $\diam(\mathcal{T}_{L-1})$-separated. 
Since 
\[
d(0,i^L(x,t_k))\leq d(0,t_k)<\diam(\mathcal{T}_{L-1}), 
\]
we can see that 
\[
i^L(x,t_1)\cap i^L(x,t_2)\cap\dots
\cap i^L(x,t_{n+1})\neq\emptyset. 
\]
Since we can write each element of $\mathcal{T}_L(x)$ 
as a union of elements in $\mathcal{T}_0(x)$, 
there exist $s_1,s_2,\dots,s_{n+1}$ in $\mathcal{T}_0(x)$ 
such that 
\[
s_1\cap s_2\cap\dots\cap s_{n+1}\neq\emptyset
\]
and $i^L(x,s_k)=i^L(x,t_k)$ for each $k=1,2,\dots,n+1$. 
We note that 
\[
\#\{i^l(x,s_k)\mid k=1,2,\dots,n+1\}=n+1
\]
for any $l>0$. 
By Theorem \ref{refine} (4), 
there exists $\zeta$ in $\mathcal{T}^{(n+1)}_0(x)$ 
such that $\zeta(k)=s_k$ for each $k=1,2,\dots,n+1$. 
This contradicts $x\notin R_\phi[B^{n+1}]$, 
and so the claim follows. 

By applying this claim to a finite set 
\[
\{\xi(k)\mid k=1,2,\dots,n\}
\cup\{\eta(k)\mid k=1,2,\dots,n\}, 
\]
we can find $L>0$ such that 
\[
\{i^L(x,\xi(k))\mid k=1,2,\dots,n\}
=\{i^L(x,\eta(k))\mid k=1,2,\dots,n\}. 
\]
It follows from Lemma \ref{j} (3) that 
there exists $\sigma\in S_d$ such that 
\[
\sigma(j_n^{L+1}(x,\xi))=j_n^{L+1}(x,\eta)\text{ and }
\sigma(\{1,2,\dots,n\})=\{1,2,\dots,n\}. 
\]
Define 
\[
N=\{j_n^{L+1}(x,\xi)(k)\mid k=1,2,\dots,d+1\}. 
\]
For $l>L$, we let 
\[
m_l=\#\{i^{l-L}(x,t)\mid t\in N\}. 
\]
It is easily seen that 
$(m_l)_{l>L}$ is a decreasing sequence of positive integers. 
Since $x$ is in $B^n$, $m_l$ is not less than $n$ for any $l>L$. 
It follows from exactly the same argument as in the claim above that 
there exists $L'>L$ such that $m_{L'}=n$. 
By the repeated use of Lemma \ref{j} (4), 
we can find $\sigma'\in S_n$ such that 
\[
\sigma'(j_n^{L'}(x,\xi))=j_n^{L'}(x,\eta). 
\]
Therefore the pair $(\phi^p(x),\phi^q(x))$ belongs to $R^n$. 
\end{proof}

\bigskip

We are now ready to give a proof of the main result. 
When $(S,\mathcal{O})$ is a topological space and 
$A$ is a subset of $S$, to simplify notation, 
we denote the induced topology on $A$ by $\mathcal{O}$, too. 

\begin{proof}[Proof of Theorem \ref{main}]
We first note that 
the equivalence relations $R^1,R^2,\dots,R^{d+1}$ and $Q^2,Q^3,\dots,Q^{d+1}$ 
are subsets of $R_\phi$. 
The proof will be completed 
by using the absorption theorem \cite[Theorem 3.2]{M3} 
and the splitting theorem \cite[Theorem 2.1]{M3} repeatedly. 

Let $\mathcal{O}_1$ be the \'etale topology on $R_\phi$. 
By Lemma \ref{R=AF} and \ref{R1minimal}, 
$(R^1,\mathcal{O}_1)$ is a minimal AF relation on $B^1=X$. 
By Lemma \ref{B2thin}, $B^2$ is a closed $R^1$-thin subset. 
By Lemma \ref{Bn+1etale}, $B^2$ is \'etale 
with respect to $(R^1,\mathcal{O}_1)$. 
By Lemma \ref{Q=AF}, 
$(Q^2,\mathcal{O}_1)$ is an AF relation on $B^2$. 
By Lemma \ref{Rn<Qn+1}, 
$R^1|B^2$ is an open subrelation of $Q^2$. 
Then \cite[Theorem 3.2]{M3} applies and yields 
an \'etale topology $\mathcal{O}_2$ on $R^1\vee Q^2$ 
satisfying the following. 
\begin{itemize}
\item[(1-a)] $(R^1\vee Q^2,\mathcal{O}_2)$ is a minimal AF relation on $X$. 
\item[(1-b)] $B^2$ is $(R^1\vee Q^2)$-thin. 
\item[(1-c)] $B^2$ is \'etale with respect to $(R^1\vee Q^2,\mathcal{O}_2)$. 
\item[(1-d)] Two topologies $\mathcal{O}_1$ and $\mathcal{O}_2$ 
agree on $Q^2$. 
\end{itemize}
By Lemma \ref{RnCn=Rphi}, one has $R^1|C^1=R_\phi|C^1$. 
In particular, $(R^1\vee Q^2)|C^1=R_\phi|C^1$. 
From Lemma \ref{RnBn+1=Rphi}, we have $R^1[B^2]=R_\phi[B^2]$. 
By Lemma \ref{RnCn=Rphi}, 
we also have $R^2|(B^2\cap C^2)=R_\phi|(B^2\cap C^2)$. 
Since $R^2$ is a subrelation of $Q^2$ by Lemma \ref{Rn<Qn}, 
one gets $Q^2|(B^2\cap C^2)=R_\phi|(B^2\cap C^2)$. 
Combining these, we have $(R^1\vee Q^2)|C^2=R_\phi|C^2$. 
Hence we obtain the following. 
\begin{itemize}
\item[(1-e)] $(R^1\vee Q^2)|(C^1\cup C^2)=R_\phi|(C^1\cup C^2)$. 
\end{itemize}

Next, we would like to apply the splitting theorem \cite[Theorem 2.1]{M3} 
to $(R^1\vee Q^2,\mathcal{O}_2)$, $B^2$ and $R^2$. 
By Lemma \ref{Rn<Qn}, 
$R^2$ is an open subrelation of $Q^2$ with the topology $\mathcal{O}_1$. 
Evidently $(R^1\vee Q^2)|B^2$ is equal to $Q^2$. 
It follows from (1-d) that 
$R^2$ is an open subrelation of $(R^1\vee Q^2)|B^2$ 
with the topology $\mathcal{O}_2$. 
Combining this with (1-a), (1-b) and (1-c), 
we can apply \cite[Theorem 2.1]{M3} and 
obtain a subrelation $\widetilde{R}^2$ of $R^1\vee Q^2$ 
satisfying the following. 
\begin{itemize}
\item[(1-f)] $(\widetilde{R}^2,\mathcal{O}_2)$ is a minimal AF relation. 
\item[(1-g)] $\widetilde{R}^2|B^2$ equals $R^2$. 
\item[(1-h)] $\widetilde{R}^2[B^2]$ equals $(R^1\vee Q^2)[B^2]$. 
\item[(1-i)] If $x\in X$ does not belong to $(R^1\vee Q^2)[B^2]$, then 
$\widetilde{R}^2[x]$ equals $(R^1\vee Q^2)[x]$. 
\item[(1-j)] Any $\widetilde{R}^2$-invariant measure on $X$ 
is $(R^1\vee Q^2)$-invariant. 
\end{itemize}
By Lemma \ref{RnCn=Rphi}, 
$R^2|(B^2\cap C^2)$ is equal to $R_\phi|(B^2\cap C^2)$. 
This together with (1-g), (1-h) and (1-i) implies that 
$\widetilde{R}^2|(C^1\cup C^2)$ is equal to $(R^1\vee Q^2)|(C^1\cup C^2)$. 
Combining this with (1-e), one gets the following. 
\begin{itemize}
\item[(1-k)] $\widetilde{R}^2|(C^1\cup C^2)=R_\phi|(C^1\cup C^2)$. 
\end{itemize}
Moreover, we claim the following. 
\begin{itemize}
\item[(1-l)] $\widetilde{R}^2[B^3]=R_\phi[B^3]$. 
\end{itemize}
To verify this, suppose that $(x,y)$ is in $R_\phi$ and $y\in B^3$. 
By $R^1[B^2]=R_\phi[B^2]$ and (1-h), 
there exists $z\in B^2$ such that $(x,z)\in\widetilde{R}^2$. 
In particular, $z$ is in $B^2\cap R_\phi[B^3]$. 
From Lemma \ref{RnBn+1=Rphi}, $z$ is in $R^2[B^3]$. 
It follows from (1-g) that 
there exists $w\in B^3$ such that $(z,w)\in\widetilde{R}^2$. 
Thus, $x$ belongs to $\widetilde{R}^2[B^3]$. 

We next wish to apply the absorption theorem \cite[Theorem 3.2]{M3} 
to $(\widetilde{R}^2,\mathcal{O}_2)$, $B^3$ and $Q^3$. 
Let us check the hypotheses. 
By (1-f), 
$(\widetilde{R}^2,\mathcal{O}_2)$ is a minimal AF relation on $X$. 
From (1-b) and (1-j), $B^3$ is $\widetilde{R}^2$-thin. 
By Lemma \ref{Bn+1etale}, 
$B^3$ is \'etale with respect to $(R^2,\mathcal{O}_1)$. 
From (1-g), $\widetilde{R}^2|B^3=R^2|B^3$. 
By (1-d), 
two topologies $\mathcal{O}_1$ and $\mathcal{O}_2$ agree on $R^2|B^3$, 
because $R^2$ is a subset of $Q^2$. 
It follows that 
$B^3$ is \'etale with respect to $(\widetilde{R}^2,\mathcal{O}_2)$. 
By Lemma \ref{Rn<Qn+1} and (1-d), 
$\widetilde{R}^2|B^3=R^2|B^3$ with the topology $\mathcal{O}_2$ 
is an open subrelation of $Q^3$ with the topology $\mathcal{O}_1$. 
Therefore we can apply \cite[Theorem 3.2]{M3} and 
obtain an \'etale topology $\mathcal{O}_3$ 
on $\widetilde{R}^2\vee Q^3$ satisfying the following. 
\begin{itemize}
\item[(2-a)] $(\widetilde{R}^2\vee Q^3,\mathcal{O}_3)$ is 
a minimal AF relation on $X$. 
\item[(2-b)] $B^3$ is $(\widetilde{R}^2\vee Q^3)$-thin. 
\item[(2-c)] $B^3$ is \'etale with respect to 
$(\widetilde{R}^2\vee Q^3,\mathcal{O}_3)$. 
\item[(2-d)] Two topologies $\mathcal{O}_1$ and $\mathcal{O}_3$ 
agree on $Q^3$. 
\end{itemize}
By Lemma \ref{RnCn=Rphi}, 
we have $R^3|(B^3\cap C^3)=R_\phi|(B^3\cap C^3)$. 
Since $R^3$ is a subrelation of $Q^3$ by Lemma \ref{Rn<Qn}, 
one gets $Q^3|(B^3\cap C^3)=R_\phi|(B^3\cap C^3)$. 
Combining this with (1-l), 
we get $(\widetilde{R}^2\vee Q^3)|C^3=R_\phi|C^3$. 
This, together with (1-k), implies the following. 
\begin{itemize}
\item[(2-e)] $(\widetilde{R}^2\vee Q^3)|(C^1\cup C^2\cup C^3)
=R_\phi|(C^1\cup C^2\cup C^3)$. 
\end{itemize}

Next, we would like to apply the splitting theorem \cite[Theorem 2.1]{M3} 
to $(\widetilde{R}^2\vee Q^3,\mathcal{O}_3)$, $B^3$ and $R^3$. 
By Lemma \ref{Rn<Qn}, 
$R^3$ is an open subrelation of $Q^3$ with the topology $\mathcal{O}_1$. 
Evidently $(\widetilde{R}^2\vee Q^3)|B^3$ is equal to $Q^3$. 
It follows from (2-d) that 
$R^3$ is an open subrelation of $(\widetilde{R}^2\vee Q^3)|B^3$ 
with the topology $\mathcal{O}_3$. 
Combining this with (2-a), (2-b) and (2-c), 
we can apply \cite[Theorem 2.1]{M3} and 
obtain a subrelation $\widetilde{R}^3$ of $\widetilde{R}^2\vee Q^3$ 
satisfying the following. 
\begin{itemize}
\item[(2-f)] $(\widetilde{R}^3,\mathcal{O}_3)$ is a minimal AF relation. 
\item[(2-g)] $\widetilde{R}^3|B^3$ equals $R^3$. 
\item[(2-h)] $\widetilde{R}^3[B^3]$ equals $(\widetilde{R}^2\vee Q^3)[B^3]$. 
\item[(2-i)] If $x\in X$ does not belong to 
$(\widetilde{R}^2\vee Q^3)[B^3]$, then 
$\widetilde{R}^3[x]$ equals $(\widetilde{R}^2\vee Q^3)[x]$. 
\item[(2-j)] Any $\widetilde{R}^3$-invariant measure on $X$ 
is $(\widetilde{R}^2\vee Q^3)$-invariant. 
\end{itemize}
By Lemma \ref{RnCn=Rphi}, 
$R^3|(B^3\cap C^3)$ is equal to $R_\phi|(B^3\cap C^3)$. 
This together with (2-g), (2-h) and (2-i) implies that 
$\widetilde{R}^3|(C^1\cup C^2\cup C^3)$ is equal to 
$(\widetilde{R}^2\vee Q^3)|(C^1\cup C^2\cup C^3)$. 
Combining this with (2-e), one gets the following. 
\begin{itemize}
\item[(2-k)] $\widetilde{R}^3|(C^1\cup C^2\cup C^3)
=R_\phi|(C^1\cup C^2\cup C^3)$. 
\end{itemize}
Moreover, we claim the following. 
\begin{itemize}
\item[(2-l)] $\widetilde{R}^3[B^4]=R_\phi[B^4]$. 
\end{itemize}
To verify this, suppose that $(x,y)$ is in $R_\phi$ and $y\in B^4$. 
By (1-l) and (2-h), 
there exists $z\in B^3$ such that $(x,z)\in\widetilde{R}^3$. 
In particular, $z$ is in $B^3\cap R_\phi[B^4]$. 
From Lemma \ref{RnBn+1=Rphi}, $z$ is in $R^3[B^4]$. 
It follows from (2-g) that 
there exists $w\in B^4$ such that $(z,w)\in\widetilde{R}^3$. 
Thus, $x$ belongs to $\widetilde{R}^3[B^4]$. 

By repeating these arguments, 
we finally obtain a subrelation $\widetilde{R}^d$ of $R_\phi$ and 
an \'etale topology $\mathcal{O}_{d+1}$ on $\widetilde{R}^d\vee Q^{d+1}$ 
satisfying the following. 
\begin{itemize}
\item[(d-a)] $(\widetilde{R}^d\vee Q^{d+1},\mathcal{O}_{d+1})$ is 
a minimal AF relation on $X$. 
\item[(d-e)] $(\widetilde{R}^d\vee Q^{d+1})|(C^1\cup C^2\cup\dots\cup C^{d+1})
=R_\phi|(C^1\cup C^2\cup\dots\cup C^{d+1})$. 
\end{itemize}
Since $C^1\cup C^2\cup\dots\cup C^{d+1}$ equals $X$, 
these two conditions imply that 
$(R_\phi,\mathcal{O}_{d+1})$ is a minimal AF equivalence relation, 
thereby completing the proof. 
\end{proof}

\appendix
\def\thesection{Appendix}
\section{A remark on free actions}
\def\thesection{\Alph{section}}

Theorem \ref{app} presented below is a generalization of a result proved 
by O. Johansen in his thesis \cite{J}, which has not been published. 
Since this is relevant for our investigation we will give a proof. 

\begin{thm}\label{app}
Let $K$ be a finite group and let $d\geq1$. 
For any free action $\phi$ of $K\oplus\Z^d$ 
as homeomorphisms on the Cantor set $X$, 
there exists a free action $\psi$ of $\Z^d$ on $X$ 
such that $R_\phi = R_\psi$. 
\end{thm}

To prove the theorem we shall need the following lemma. 

\begin{lem}
Let $K$ be a finite group acting freely on the Cantor set $X$. 
There exists a clopen subset $E$ in $X$ such that 
$\{k(E)\mid k\in K\}$ is a (clopen) partition of $X$. 
\end{lem}
\begin{proof}
For each $x$ in $X$, 
a clopen set $E_x$ containing $x$ can be chosen 
such that the sets $k(E_x)$, for $k$ in $K$, are disjoint. 
Then $\{E_x\mid x\in X\}$ covers $X$. 
Choose a finite subcovering $\{E_{x_j}\mid j\in J\}$. 
Let $\mathcal{P}$ be the clopen partition of $X$ 
generated by $\{k(E_{x_j})\mid k\in K, \ j\in J\}$. 
Then $\mathcal{P}$ is $K$-invariant, i.e. 
$k(M)\in\mathcal{P}$ if $M\in\mathcal{P}$ and $k\in K$. 
Furthermore, if $M\in\mathcal{P}$ 
then $M\subset E_{x_j}$ for some $j\in J$, 
and so the sets $k(M), k\in K$, are disjoint. 
The relation $M\sim N$ on $\mathcal{P}$ 
defined by $M=k(N)$ for some $k\in K$ is an equivalence relation. 
Choose one element $M_l$ from each equivalence class and 
let $E=\bigcup_lM_l$. 
Then $E$ will have the desired property. 
\end{proof}

\begin{proof}[Proof of Theorem \ref{app}]
It will be convenient to denote 
$\phi^{(k,0)}$ and $\phi^{(0,h)}$ by $k$ and $h$, respectively, 
i.e. $\phi^{(k,0)}(x)=k(x)$ and $\phi^{(0,h)}(x)=h(x)$, 
where $k\in K$, $h\in\Z^d$ and $x\in X$. 
Let $Y\subset X$ be the clopen set of the lemma 
associated to the action of $K$. 
Let $K=\{e{=}k_0,k_1,\dots,k_{n-1}\}$, 
where $e$ is the identity element of $K$, 
and let $\{Y_i{=}k_i(Y)\mid i=0,1,\dots,n-1\}$ be 
the clopen partition of $X$ according to the lemma. 

Let $\pi_Y:X\to Y$ be the continuous map defined by 
\[
\pi_Y(x)=k^{-1}_i(x) \qquad \text{if }x\in Y_i, \quad i=0,1,\dots,n-1. 
\]
It is easy to see that 
$\pi_Y\circ h\circ k$ equals $\pi_Y\circ h$ for any $h\in\Z^d$ and $k\in K$. 
Let $\{h_1,h_2,\dots,h_d\}$ be a basis for $\Z^d$. 
For $l\in\{2,\dots,d\}$, let $\psi_l:X\to X$ be defined by 
\[
\psi_l(x)=k_j(\pi_Y(h_l(x))) \qquad\text{if }x\in Y_j 
\]
and let 
\[
\psi_1(x)=
\begin{cases}
\pi_Y(h_1(x)) & \text{if }x\in Y_{n-1} \\
k_{i+1}(k_i^{-1}(x)) & \text{if }x\in Y_i\text{ and }i\neq n-1. 
\end{cases}
\]
Clearly $\psi_1,\psi_2,\dots,\psi_d$ are homeomorphisms on $X$. 

Now let $2\leq l_1,l_2\leq d$. 
For $x\in Y_i$, we have 
\begin{align*}
\psi_{l_1}(\psi_{l_2}(x))
&=(\psi_{l_1}\circ k_i\circ\pi_Y\circ h_{l_2})(x) \\
&=(k_i\circ\pi_Y\circ h_{l_1}\circ k_i\circ\pi_Y\circ h_{l_2})(x) \\
&=(k_i\circ\pi_Y\circ h_{l_1}\circ h_{l_2})(x) \\
&=(k_i\circ\pi_Y\circ h_{l_2}\circ h_{l_1})(x)=\psi_{l_2}(\psi_{l_1}(x)). 
\end{align*}
In a similar way one checks that 
$\psi_1$ commutes with $\psi_l$, where $l\in\{2,\dots,d\}$. 
So $\psi_1,\psi_2,\dots,\psi_d$ give rise to 
an action $\psi$ of $\Z^d$ as homeomorphisms on $X$. 

We prove that $R_\phi=R_\psi$, or, 
what is the same, that the $\phi$-orbits and the $\psi$-orbits coincide. 
Obviously we have $R_\psi[x]\subset R_\phi[x]$ for all $x\in X$. 
To prove the converse, let $x\in Y_i$. 
Take $k\in K$ and suppose $k(x)\in Y_j$.  
Then one can see $\psi_1^{j-i}(x)=k(x)$, 
and so $k(x)$ is in $R_\psi[x]$. 
Next, to prove $h_l(x)\in R_\psi[x]$, 
assume that $h_l(x)$ belongs to $Y_j$. 
Then, for $l=2,\dots,d$, 
\[
h_l(x)=(k_j\circ\pi_Y\circ h_l)(x)
=(k_j\circ k_i^{-1}\circ\psi_l)(x), 
\]
and 
\[
h_1(x)=(k_j\circ\pi_Y\circ h_1\circ \psi_1^{n-1-i})(x)
=(k_j\circ\psi_1^{n-i})(x). 
\]
Therefore $h_l(x)\in R_\psi[x]$. 

What remains to be shown is that 
$\psi$ is a free action of $\Z^d$. 
This is easy to see. 
In fact, if $x\in Y_j$, then 
each $\psi_2,\dots,\psi_d$ map $x$ into $Y_j$, 
and $\psi_1$ maps $x$ into $Y_{j+1}$, 
where $j+1$ is understood modulo $n$. 
Suppose that 
$(\psi_1^{a_1}\circ\psi_2^{a_2}\circ\dots\circ\psi_d^{a_d})(x)=x$. 
Then $a_1$ is divisible by $n$. 
Let $a_1=nb$. 
Since $K$ commutes with every $h_1,h_2,\dots,h_d$, 
it is easily seen that 
there exists $k\in K$ such that 
\[
x=(\psi_1^{a_1}\circ\psi_2^{a_2}\circ\dots\circ\psi_d^{a_d})(x)
=(k\circ h_1^b\circ h_2^{a_2}\circ\dots\circ h_d^{a_d})(x). 
\]
It follow from the freeness of $\phi$ that 
$b=a_2=\dots=a_d=0$. 
Hence $\psi$ is free. 
\end{proof}

The following corollary follows 
from Theorem \ref{main} and Theorem \ref{app}.

\begin{cor}
Any minimal action of a finitely generated
abelian group on the Cantor set is affable. 
\end{cor}


\begin{thebibliography}{ABCD}
\bibitem[BBG]{BBG}
J. Bellissard, R. Benedetti and J. -M. Gambaudo, 
\textit{Spaces of tilings, finite telescopic approximations 
and gap-labeling}, 
Comm. Math. Phys. 261 (2006), 1--41. 
\bibitem[BG]{BG}
R. Benedetti and J. -M. Gambaudo, 
\textit{On the dynamics of $G$-solenoids. Applications to Delone sets}, 
Ergodic Theory Dynam. Systems 23 (2003), 673--691. 
\bibitem[CDEFT]{CDEFT}
Siu-Wing Cheng, Tamal K. Dey, Herbert Edelsbrunner, 
Michael A. Facello and Shang-Hua Teng, 
\textit{Sliver exudation}, 
J. ACM 47 (2000), 883--904. 
\bibitem[CFW]{CFW}
A. Connes, J. Feldman and B. Weiss, 
\textit{An amenable equivalence relation is generated 
by a single transformation}, 
Ergodic Theory Dynamical Systems 1 (1981), 431--450. 
\bibitem[DJK]{DJK}
R. Dougherty, S. Jackson and A. S. Kechris, 
\textit{The structure of hyperfinite Borel equivalence relations}, 
Trans. Amer. Math. Soc. 341 (1994), 193--225. 
\bibitem[D]{D}
H. A. Dye, 
\textit{On groups of measure preserving transformation. I}, 
Amer. J. Math. 81 (1959), 119--159. 
\bibitem[F]{F}
A. Forrest, 
\textit{A Bratteli diagram for commuting homeomorphisms of the Cantor set}, 
Internat. J. Math. 11 (2000), 177--200. 
\bibitem[GJ]{GJ}
S. Gao and S. Jackson, 
\textit{Countable abelian group actions and 
hyperfinite equivalence relations}, 
preprint. 
http://www.math.unt.edu/~sgao/pub/paper30.html
\bibitem[GMPS1]{GMPS1}
T. Giordano, H. Matui, I. F. Putnam and C. F. Skau, 
\textit{Orbit equivalence for Cantor minimal $\Z^2$-systems}, 
J. Amer. Math. Soc. 21 (2008), 863--892. 
math.DS/0609668
\bibitem[GMPS2]{GMPS2}
T. Giordano, H. Matui, I. F. Putnam and C. F. Skau, 
\textit{The absorption theorem for affable equivalence relations}, 
Ergodic Theory Dynam. Systems 28 (2008), 1509--1531. 
arXiv:0705.3270
\bibitem[GPS1]{GPS1}
T. Giordano, I. F. Putnam and C. F. Skau, 
\textit{Topological orbit equivalence and $C^*$-crossed products}, 
J. Reine  Angew. Math. 469 (1995), 51--111. 
\bibitem[GPS2]{GPS2}
T. Giordano, I. F. Putnam and C. F. Skau, 
\textit{Affable equivalence relations and orbit structure 
of Cantor dynamical systems}, 
Ergodic Theory Dynam. Systems 24 (2004), 441--475. 
\bibitem[GPS3]{GPS3}
T. Giordano, I. F. Putnam and C. F. Skau, 
\textit{Cocycles for Cantor minimal $\Z^d$-systems}, 
to appear in Internat. J. Math. 
\bibitem[HPS]{HPS}
R. H. Herman, I. F. Putnam and C. F. Skau, 
\textit{Ordered Bratteli diagrams, dimension groups and 
topological dynamics}, 
Internat. J. Math. 3 (1992), 827--864. 
\bibitem[JKL]{JKL}
S. Jackson, A. S. Kechris and A. Louveau, 
\textit{Countable Borel equivalence relations}, 
J. Math. Log. 2 (2002), 1--80. 
\bibitem[J]{J}
\O. Johansen, 
\textit{Ordered $K$-theory and Bratteli diagrams: 
Implications for Cantor minimal systems}, 
Ph.D. thesis, NTNU, 1998. 
\bibitem[KM]{KM}
A. S. Kechris and B. D. Miller, 
\textit{Topics in Orbit Equivalence}, 
Lecture Notes in Math. 1852, Springer-Verlag, 
Berlin Heidelberg New York, 2004. 
\bibitem[KP]{KP}
J. Kellendonk and I. F. Putnam, 
\textit{Tilings, $C^*$-algebras, and $K$-theory}, 
Directions in Mathematical Quasicrystals, 177--206, 
CRM Monogr. Ser. 13, Amer. Math. Soc., Providence, RI, 2000. 
\bibitem[K]{K}
W. Krieger, 
\textit{On ergodic flows and the isomorphism of factors}, 
Math. Ann. 223 (1976), 19--70. 
\bibitem[Mac]{Mac}
G. W. Mackey, 
\textit{Ergodic theory and virtual groups}, 
Math. Ann. 166 (1966) 187--207. 
\bibitem[Mat]{Mat}
J. Matou\u{s}ek, 
\textit{Lectures on discrete geometry}, 
Graduate Texts in Mathematics, 212. Springer-Verlag, New York, 2002. 
\bibitem[M1]{M1}
H. Matui, 
\textit{Affability of equivalence relations 
arising from two-dimensional substitution tilings}, 
Ergodic Theory Dynam. Systems 26 (2006), 467--480. 
math.DS/0506251
\bibitem[M2]{M2}
H. Matui, 
\textit{A short proof of affability 
for certain Cantor minimal $\Z^2$-systems}, 
Canad. Math. Bull. 50 (2007), 418--426. 
math.DS/0506250
\bibitem[M3]{M3}
H. Matui, 
\textit{An absorption theorem 
for minimal AF equivalence relations on Cantor sets}, 
J. Math. Soc. Japan 60 (2008), 1171--1185. 
arXiv:0712.0733
\bibitem[OW]{OW}
D. S. Ornstein and B. Weiss, 
\textit{Ergodic theory of amenable group actions I: 
The Rohlin lemma}, 
Bull. Amer. Math. Soc. 2 (1980), 161--164. 
\bibitem[P]{P}
I. F. Putnam, 
\textit{Orbit equivalence of Cantor minimal systems: 
a survey and a new proof}, 
preprint. 
\bibitem[R]{R}
J. Renault, 
\textit{A Groupoid Approach to $C^*$-algebras}, 
Lecture Notes in Mathematics 793, Springer, Berlin, 1980. 
\bibitem[Sch]{Sch}
M. Schlottmann, 
\textit{Periodic and quasi-periodic Laguerre tilings}, 
Internat. J. Modern Phys. B 7 (1993), 1351--1363. 
\bibitem[St]{St}
R. P. Stanley, 
\textit{Decompositions of rational convex polytopes}, 
Ann. Discrete Math. 6 (1980), 333--342. 
\end{thebibliography}
\end{document}